\documentclass[aos,MSNbibl,nameyear,seceqn,dvips]{arximspdf}
\usepackage{dcolumn}
\usepackage{graphicx}

\doi{10.1214/14-AOS1210} 
\volume{42}
\issue{3}
\pubyear{2014}
\firstpage{970}
\lastpage{1002}

\makeatletter
\newcolumntype{d}[1]{D{.}{.}{#1}}
\newcommand{\widebar}{\overline}
\newcommand{\rrVert}{\Vert}
\newcommand{\rrvert}{\vert}
\newcommand{\llVert}{\Vert}
\newcommand{\llvert}{\vert}
\newproclaim{rem}{Remark}
\newtheorem{teo}{Theorem}
\newtheorem{pro}{Proposition}
\newtheorem{lemma}{Lemma}
\newtheorem{coro}{Corollary}

\makeatother

\begin{document}
\begin{frontmatter}

\title{Nonparametric maximum likelihood approach to multiple
change-point problems\thanksref{T1}}
\runtitle{Nonparametric detection of multiple change-points}

\begin{aug}
\author[A]{\fnms{Changliang} \snm{Zou}\ead[label=e1]{nk.chlzou@gmail.com}},
\author[B]{\fnms{Guosheng} \snm{Yin}\corref{}\ead[label=e2]{gyin@hku.hk}},
\author[C]{\fnms{Long} \snm{Feng}\ead[label=e3]{flnankai@126.com}}\\
\and
\author[A]{\fnms{Zhaojun} \snm{Wang}\ead[label=e4]{zjwang@nankai.edu.cn}}
\runauthor{Zou, Yin, Feng and Wang}
\affiliation{Nankai University, University of Hong Kong, Nankai
University\\ and Nankai University}
\address[A]{C. Zou\\
Z. Wang\\
Institute of Statistics\\
Nankai University\\
Tianjin 300071\\
China\\
\printead{e1}\\
\phantom{E-mail:\ }\printead*{e4}}
\address[B]{G. Yin\\
Department of Statistics\\
\quad and Actuarial Science\\
University of Hong Kong\\
Hong Kong\\
\printead{e2}}
\address[C]{L. Feng\\
School of Mathematical Sciences\\
Nankai University\\
Tianjin 300071\\
China\\
\printead{e3}}
\end{aug}
\thankstext{T1}{Supported in part by the NNSF of China Grants 11131002, 11101306, 11371202,
the RFDP of China Grant 20110031110002, the Foundation for the Author of National Excellent Doctoral
Dissertation of PR China 201232, New Century Excellent Talents in University and
also by a Grant (784010) from the Research grants Council of Hong Kong.}

\received{\smonth{11} \syear{2013}}
\revised{\smonth{2} \syear{2014}}

%
\begin{abstract}
In multiple change-point problems, different data segments often follow
different distributions, for which the changes may occur in the mean,
scale or
the entire distribution from one segment to another. Without the
need to know the number of change-points in advance, we propose a
nonparametric maximum likelihood approach to detecting multiple
change-points. Our method does not impose any parametric assumption
on the underlying distributions of the data sequence, which is thus
suitable for detection of any changes in the distributions. The
number of change-points is determined by the Bayesian information
criterion and the locations of the change-points can be estimated via
the dynamic programming algorithm and
the use of the intrinsic order structure of the likelihood function. Under
some mild conditions, we show that the new method provides
consistent estimation with an optimal rate. We also suggest a
prescreening procedure to
exclude most of the irrelevant points prior to the implementation of
the nonparametric likelihood method. Simulation studies show
that the proposed method has satisfactory performance of identifying
multiple change-points in terms of estimation accuracy and
computation time.
\end{abstract}

%

\end{frontmatter}

\section{Introduction}\label{sec1}
The literature devoted to change-point models is vast, particularly
in the areas of economics, genome research, quality control, and
signal processing. When there are notable changes in a sequence of
data, we can typically break the sequence into several data
segments, so that the observations within each segment are
relatively homogeneous. In the conventional change-point problems,
the posited models for different data segments are often of the same
structure but with different parameter values.
However, the underlying distributions are typically unknown, and thus
parametric methods potentially suffer from model misspecification.
The least-squares fitting is the standard choice for the MCP, while
its performance often deteriorates when the error follows a
heavy-tailed distribution or when the data contain outliers.

Without imposing any parametric modeling assumption, we consider the multiple
change-point problem (MCP) based on independent data
$\{X_i\}_{i=1}^n$, such that
%
%
\begin{equation}
\label{om1} \qquad X_{i}\sim F_k(x), \qquad\tau_{k-1}
\leq i \leq\tau_k-1,  k=1,\ldots,K_n+1; i=1,\ldots,n,
\end{equation}
where $K_n$ is the true number of change-points, $\tau_k$'s are the
locations of these change-points with the convention of $\tau_0=1$
and $\tau_{K_n+1}=n+1$, and $F_k$ is the cumulative distribution function (C.D.F.) of segment $k$ satisfying $F_k\neq F_{k+1}$.
The number of change-points $K_n$ is allowed to grow with the sample
size $n$.

Although extensive research has been conducted to estimate the
number of change-points $K_n$ and the locations of these
change-points $\tau_k$'s, most of the work assumes that $F_k$'s
belong to some-known parametric functional families or that they
differ only in their locations (or scales). For a comprehensive
coverage on single change-point problems $(K_n=1)$, see \citet{CsoHor97}. The standard
approach to the MCP is based on least-squares or likelihood methods
via a dynamic programming (DP) algorithm in conjunction with a
selection procedure such as the Bayesian information criterion (BIC)
for determining the number of change-points
[\citet{Yao88}; \citet{YaoAu89}; \citet{CheGup97}; \citeauthor{BaiPer98} (\citeyear{BaiPer98,BaiPer03});
\citet{BraBraMul00}; \citet{Haw01}; \citet{Lav05}]. By reframing the
MCP in a
variable selection context, \citet{HarLev10} proposed a penalized
least-squares criterion with a
LASSO-type penalty [\citet{Tib96}]. \citet{CheZha12}
developed a graph-based approach to detecting change-points, which
is applicable in high-dimensional data and non-Euclidean data. Other
recent development in this area includes \citet{Rig10}, \citet{KilFeaEck12} and \citet{ArlCelHar12}.

Our goal is to develop an efficient nonparametric procedure for the MCP in
(\ref{om1}) without imposing any parametric structure on the
$F_k$'s; virtually any salient difference between two successive
C.D.F.'s (say, $F_k$ and $F_{k+1}$) would ensure detection of the
change-point asymptotically. In the nonparametric context, most of the
existing work focuses on the single change-point problem by using
some seminorm on the difference between pre- and post-empirical
distributions at the change-point [\citet{Dar76}; \citet{Car88};
\citet{Dum91}]. \citet{Gua04} studied a semiparametric
change-point model based on the empirical likelihood, and applied
the method to detect the change from a distribution to a weighted
one. \citet{Zouetal07} proposed another empirical likelihood
approach without assuming any relationship between the two
distributions. However, extending these methods to the MCP is not
straightforward. \citet{Lee96} proposed to use the weighted empirical
measure to detect two different nonparametric distributions over a
window of observations and then run the window through the full data sequence
to detect the number of change-points. Although the approach of
\citet{Lee96}
is simple and easy to implement,
our simulation studies show that
even with elaborately chosen tuning parameters the estimates of the
locations $\tau_k$'s as well as the number of change-points are not
satisfactory. This may be partly due to the ``local'' nature of\vadjust{\goodbreak}
the running window, and thus the information in the data is not
fully and efficiently utilized. \citet{MatJam}
proposed a new estimation method, ECP, under multivariate settings,
which is based on hierarchical clustering by recursively using a
single change-point estimation procedure.

Observing the connection between multiple change-points and
goodness-of-fit tests, we propose a nonparametric maximum likelihood
approach to the MCP. Our proposed nonparametric multiple
change-point detection (NMCD) procedure can be regarded as a
nonparametric counterpart of the classical least-squares MCP method
[\citet{Yao88}]. Under some mild conditions, we demonstrate that the NMCD can
achieve the optimal rate, $O_p(1)$, for the estimation of the
change-points without any distributional assumptions. Due to the use
of empirical \mbox{distribution} functions, technical arguments for
controlling the supremum of the nonparametric likelihood function
are nontrivial and are interesting in their own rights. As a
matter of fact, some techniques regarding the empirical process
have been nicely integrated with
the MCP methodologies. In addition, our theoretical results are applicable
to the situation with a diverging number of change-points, that is,
when the
number of change-points, $K_n$, grows as $n$ goes to infinity. This
substantially enlarges the scope of applicability of the proposed
method, from a traditional fixed dimensionality to a more challenging
high-dimensional setting.

In the proposed NMCD procedure,
the number of change-points, $K_n$, is determined by the BIC.
Given $K_n$, the DP algorithm utilizes
the intrinsic order structure of the likelihood to recursively
compute the maximizer of the objective function with a complexity of
$O(K_nn^2)$. To exclude most of the irrelevant points, we also suggest
an initial screening procedure so that the
NMCD is implemented in a much lower-dimensional space.
Compared with existing parametric and nonparametric approaches,
the proposed NMCD has satisfactory performance of
identifying multiple change-points in terms of estimation accuracy and
computation time. It offers robust and effective detection capability
regardless of whether the $F_k$'s differ in the location, scale, or shape.

The remainder of the paper is organized as follows. In Section~\ref
{sec2}, we
first describe how to recast the MCP in (\ref{om1}) into a
maximization problem and then introduce our nonparametric likelihood
method followed by its asymptotic properties. The algorithm and
practical implementation are presented in Section~\ref{sec3}. The numerical
performance and comparisons with other existing methods are
presented in Section~\ref{sec4}. Section~\ref{sec5} contains a real
data example to
illustrate the application of our NMCD method. Several remarks draw
the paper to its conclusion in Section~\ref{sec6}. Technical proofs are
provided in the \hyperref[app]{Appendix}, and the proof of a corollary and additional
simulation results are given in the supplementary material [\citet{Zouetal}].

\section{Nonparametric multiple change-point detection}\label{sec2}
\subsection{NMCD method}\label{sec2.1}
Assume that $Z_1,\ldots,Z_n$ are independent and identically
distributed from $F_0$, and let ${\widehat{F}}_n$ denote the empirical
C.D.F. of the
sample, then $n{\widehat{F}}_n(u)\sim\operatorname{Binomial}(n,F_0(u))$. If
we regard
the sample as binary data with the probability of success ${\widehat{F}}_n(u)$,
this leads to the nonparametric maximum log-likelihood
\[
n \bigl\{{\widehat{F}}_n(u)\log\bigl({{\widehat{F}}_n(u)}
\bigr)+\bigl(1-{\widehat{F}}_n(u)\bigr)\log\bigl({1-{
\widehat{F}}_n(u)} \bigr) \bigr\}.
\]
In the context of (\ref{om1}), we can write the joint
log-likelihood for a candidate set of change-points
$(\tau_1'<\cdots<\tau_L')$ as
%
%
\begin{eqnarray}
\label{sjl} {\mathcal L}_u\bigl(\tau_1',\ldots,\tau_L'\bigr)&=&\sum_{k=0}^L
\bigl(\tau_{k+1}'-\tau_k'\bigr)
\bigl\{{\widehat{F}}_{\tau_k'}^{\tau_{k+1}'}(u) \log\bigl({{
\widehat{F}}_{\tau_k'}^{\tau_{k+1}'}(u)} \bigr)
\nonumber\\[-8pt]\\[-8pt]
&&{}+ \bigl(1-{\widehat{F}}_{\tau_k'}^{\tau_{k+1}'}(u) \bigr) \log\bigl({1-{
\widehat{F}}_{\tau_k'}^{\tau_{k+1}'}(u)} \bigr) \bigr\},\nonumber
\end{eqnarray}
where ${\widehat{F}}_{\tau_{k}'}^{\tau_{k+1}'}(u)$ is the empirical
C.D.F. of
the subsample $\{X_{\tau_{k}'},\ldots,X_{{\tau_{k+1}'}-1}\}$ with
$\tau_0'=1$ and $\tau_{L+1}'=n+1$.
To estimate the change-points $1<\tau_1'<\cdots<\tau_L'\leq n$, we can
maximize (\ref{sjl}) in an integrated form
%
%
\begin{equation}
\label{pl} R_n\bigl(\tau_1',\ldots,
\tau_L'\bigr)= \int_{-\infty}^{\infty}{
\mathcal L}_u\bigl(\tau_1',\ldots,
\tau_L'\bigr)\,dw(u),
\end{equation}
where $w(\cdot)$ is some positive weight function so that
$R_n(\cdot)$ is finite, and the integral is used to combine all the
information across $u$. The rationale of using (\ref{pl}) can be
clearly seen from the behavior of its population counterpart.
For simplicity, we assume that there exists only one change-point
$\tau_1$, and let $\tau_1/n\rightarrow q_1\in(0,1)$ and $\tau
_1'/n\rightarrow
\theta\in(0,1)$. Through differentiation with respect to $\theta$,
it can be verified that the limiting function of ${\mathcal L}_u(\tau_1')/n$,
\begin{eqnarray*}
Q_{u}(\theta)&=&\theta \bigl\{F_{\theta}^{(1)}(u)\log
\bigl(F_{\theta}^{(1)}(u) \bigr)+ \bigl(1-F_{\theta}^{(1)}(u)
\bigr)\log\bigl({1-F_{\theta}^{(1)}(u)} \bigr) \bigr\}
\\
&&{}+(1-\theta) \bigl\{F_{\theta}^{(2)}(u)\log\bigl(F_{\theta
}^{(2)}(u)
\bigr)+ \bigl(1-F_{\theta}^{(2)}(u)\bigr)\log
\bigl({1-F_{\theta}^{(2)}(u)} \bigr) \bigr\},
\end{eqnarray*}
increases as $\theta$ approaches $q_1$ from both sides, where
\begin{eqnarray*}
F_{\theta}^{(1)}(u)&=&\frac{\min(q_1,\theta)F_1(u)+\max(\theta
-q_1,0)F_{2}(u)} {
\min(q_1,\theta)+\max(\theta-q_1,0)}\quad\mbox{and}
\\
F_{\theta}^{(2)}(u)&=&\frac{\max(q_1-\theta,0)F_1(u)+\min(1-\theta,1-q_1)F_{2}(u)} {
\max(q_1-\theta,0)+\min(1-\theta,1-q_1)},
\end{eqnarray*}
are\vspace*{-2pt} the limits of ${\widehat{F}}_{1}^{\tau_{1}'}(u)$ and
${\widehat{F}}_{\tau_1'}^{n+1}(u)$, respectively.
This implies that the function
$\int_{-\infty}^{\infty}Q_{u}(\theta)\,dw(u)$ attains its local
maximum at the true location of the change-point, $q_1$.

\begin{rem}\label{rem1}
The log-likelihood function (\ref{sjl}) is
essentially
related to the two-sample goodness-of-fit (GOF) test statistic
based on the nonparametric likelihood ratio [\citet{EinMcK03}; \citet{Zha06}]. To
see this, let $Z_1,\ldots,Z_n$ be
independent, and suppose that $Z_1,\ldots, Z_{n_1}$ have a common
continuous distribution function $F_1$, and $Z_{n_1+1},\ldots,Z_n$ have
$F_2$. We are interested in testing the null hypothesis $H_0$ that
$F_1(u)=F_2(u)$ for all $u\in(-\infty,\infty)$ against $H_1$ that
$F_1(u)\neq F_2(u)$ for some $u\in(-\infty,\infty)$.
For each fixed $u\in(-\infty,\infty)$, a natural approach is to
apply the likelihood ratio test,
\begin{eqnarray*}
G_u&=& n_1 \biggl\{{\widehat{F}}_{1}^{n_1+1}(u)
\log\biggl(\frac{{{\widehat{F}}
_{1}^{n_1+1}(u)}}{{\widehat{F}}_n(u)} \biggr)+ \bigl(1-{\widehat
{F}}_{1}^{n_1+1}(u)
\bigr)\log\biggl(\frac{1-{\widehat
{F}}_{1}^{n_1+1}(u)}{1-{{\widehat{F}}
_n(u)}} \biggr) \biggr\}
\\
&&{}+n_2 \biggl\{{\widehat{F}}_{n_1+1}^{n+1}(u)
\log\biggl(\frac{{\widehat{F}}
_{n_1+1}^{n+1}(u)}{{{\widehat{F}}_n(u)}} \biggr)+ \bigl(1-{\widehat
{F}}_{n_1+1}^{n+1}(u)
\bigr)\log\biggl(\frac{1-{\widehat{F}}
_{n_1+1}^{n+1}(u)}{1-{{\widehat{F}}_n(u)}} \biggr) \biggr\},
\end{eqnarray*}
where ${\widehat{F}}_n(u)$ corresponds to the C.D.F. of the pooled
sample. By
noting that $n_1{\widehat{F}}_{1}^{n_1+1}(u)+n_2{\widehat
{F}}_{n_1+1}^{n+1}(u)=n{\widehat{F}}_n(u)$,
$G_u$ would be of the same form as (\ref{sjl}) with $L=1$ up to a
constant which does not depend on the segmentation point~$n_1$.
\citet{EinMcK03} considered using $G_u$ to test whether
there is at most one change-point.

In the two-sample GOF test, \citeauthor{Zha02} (\citeyear{Zha02,Zha06})
demonstrated that by choosing appropriate weight functions $w(u)$ we
can produce new omnibus tests that are generally much more powerful
than the conventional ones such as Kolmogorov--Smirnov,
Cram\'{e}r--von Mises and Anderson--Darling test statistics.
If we take $dw(u)=\{{\widehat{F}}_n(u)(1-{\widehat{F}}_n(u))\}
^{-1}\,d{\widehat{F}}_n(u)$,
and also note that ${\mathcal L}_u$ is zero for $u\in(-\infty,X_{(1)})$
and $u\in(X_{(n)},\infty)$ where $X_{(1)}< \cdots< X_{(n)}$
represent the order statistics, the objective function in (\ref{pl})
can be rewritten as
%
%
\begin{eqnarray}\label{fpl}
&& R_n\bigl(\tau_1',\ldots,
\tau_L'\bigr)\nonumber
\\
&&\qquad = \int_{X_{(1)}}^{X_{(n)}}{
\mathcal L}_u\bigl(\tau_1',\ldots,
\tau_L'\bigr)\bigl\{{\widehat{F}}_n(u)
\bigl(1-{\widehat{F}}_n(u)\bigr)\bigr\}^{-1}\,d{\widehat{F}}
_n(u)
\\
&&\qquad = n\sum_{k=0}^L\sum
_{l=2}^{n-1}\bigl(\tau_{k+1}'-
\tau_k'\bigr)\frac
{{\widehat{F}}
_{kl}\log
{\widehat{F}}_{kl}+(1-{\widehat{F}}_{kl})\log(1-{\widehat
{F}}_{kl})}{l(n-l)},\nonumber
\end{eqnarray}
where\vspace*{-1pt} ${\widehat{F}}_{kl}={\widehat{F}}_{\tau_k'}^{\tau
_{k+1}'}(X_{(l)})$. As recommended
by \citet{Zha02}, we take a common ``continuity correction'' by
replacing ${\widehat{F}}_{kl}$ with ${\widehat{F}}_{kl}-1/\{2(\tau
_{k+1}'-\tau_k')\}$ for all
$k$~and~$l$.

To determine $L$ in the MCP, we observe
that $Q_u({\theta})$ is a convex function with respect to
$\theta$, and thus
\[
\max_{\tau_1'<\cdots<\tau_L'} R_n\bigl(\tau_1',\ldots,\tau_L'\bigr)\leq
\max_{\tau_1'<\cdots<\tau_{L+1}'}
R_n\bigl(\tau_1',\ldots,
\tau_{L+1}'\bigr),
\]
which\vspace*{2pt} means that the maximum
log-likelihood $\max_{\tau_1'<\cdots<\tau_L'}
R_n(\tau_1',\ldots,\tau_L')$ is a nondecreasing function in $L$.
Hence, we can use Schwarz's Bayesian information criterion (BIC) to
strike a balance between the likelihood and the number of
change-points by incorporating a penalty for large $L$.
More specifically, we identify the value of $L$
by minimizing
%
%
\begin{equation}
\label{BIC} \mathrm{BIC}_{L}=-\max_{\tau_1'<\cdots<\tau_L'}
R_n\bigl(\tau_1',\ldots,
\tau_L'\bigr)+L\zeta_n
\end{equation}
and $\zeta_n$ is a proper sequence going to infinity. \citet{Yao88}
used the BIC with $\zeta_n=\log n$ to select the number of
change-points and showed its consistency in the least-squares
framework. However, the traditional BIC
tends to select a model with some
spurious change-points. Detailed discussions on the choice of
$\zeta_n$ and other tuning parameters are given in Section~\ref{sec3.2}.
\end{rem}

\subsection{Asymptotic theory}\label{sec2.2}

In the context of change-point estimation, it is well known
that the points around the true change-point cannot be distinguished
asymptotically with a fixed change magnitude. In the least-squares
fitting, the total \mbox{variation} with perfect segmentation is
asymptotically equivalent to that with an estimate of the
change-point in a neighborhood of the true change-point [\citet{YaoAu89}].
For example, suppose that there is only one change-point
$\tau$ with a change size $\delta$, then we can only achieve
$\delta^2|\hat{\tau}_{\mathrm{MLE}}-\tau|=O_p(1)$ as $n\rightarrow
\infty$,
where $\hat{\tau}_{\mathrm{MLE}}$ denotes the maximum likelihood
estimator (MLE) of $\tau$ [see Chapter~1 of \citet{CsoHor97}]. For single
change-point nonparametric models, \citet{Dar76} obtained a rate
of $o_p(n)$, \citet{Car88} derived a rate of $O(n^{\alpha})$ a.s.
(almost surely) for any $\alpha>1/2$, and \citet{Dum91} achieved a rate of $O_p(1)$. The estimator in \citet{Lee96} is
shown to be consistent a.s.
and the differences between the estimated and true locations of
change-points are of order $O(\log n)$ a.s.

Let $\mathcal{G}_n(L)=\{\hat{\tau}_1,\ldots,\hat{\tau}_{L}\}$ denote
the set of estimates of the change-points using the proposed NMCD. The next
theorem establishes the desirable property for the NMCD estimator
when $K_n$ is prespecified---$\mathcal{G}_n(K_n)$ is asymptotically close
to the true change-point set. Let $C_{K_n}(\delta_n)$ contain
all the sets in the $\delta_n$-neighborhood of the true locations,
\begin{eqnarray*}
&&C_{K_n}(\delta_n)
\\
&&\qquad = \bigl\{\bigl(\tau_1',\ldots,\tau_{{K_n}}'
\bigr)\dvtx 1<\tau_1'<\cdots<\tau_{K_n}'
\leq n, \bigl|\tau_s'-\tau_s\bigr|\leq\delta
_n \mbox{ for } 1\leq s\leq K_n \bigr\},
\end{eqnarray*}
where $\delta_n$ is some positive sequence. Denote
$F_{k,\theta}=\theta F_k+(1-\theta)F_{k+1}$ for $0<\theta<1$. For
$r=1,\ldots, K_n$, define
\[
\eta(u;F_r,F_{r,\theta})= F_{r}(u)\log\biggl(
\frac{F_{r}(u)}{F_{r,\theta}(u)} \biggr)+ \bigl(1-F_{r}(u)\bigr)\log
\biggl(
\frac{1-F_{{r}}(u)}{1-F_{{r},\theta
}(u)} \biggr),
\]
which is the Kullback--Leibler distance between two Bernoulli
distributions with respective success probabilities $F_{r}(u)$ and
$F_{{r,\theta}}(u)$. Hence, whenever $F_{r}(u)\neq F_{{r+1}}(u)$,
and accordingly $F_{r}(u)\neq F_{r,\theta}(u)$,
$\eta(u;F_r,F_{r,\theta})$ is strictly larger than zero.
Furthermore, for $r=1,\ldots, K_n$, define
\[
\eta_r(u)=\eta(u;F_r,F_{r,1/2})+
\eta(u;F_{r+1},F_{r,1/2}).
\]

To establish the consistency of the proposed NMCD, the following
assumptions are imposed:
\begin{longlist}[(A2)]
\item[(A1)] $F_1,\ldots, F_{K_n+1}$ are continuous and $F_k\neq
F_{k+1}$ for $k=1,\ldots,K_n$.

\item[(A2)]
Let $\lambda_n=\min_{1\le k \le K_n+1} (\tau_{k}-\tau_{k-1})$;
$\lambda_n \to\infty$ as $n\to\infty$.

\item[(A3)] ${\widehat{F}}_n(u)\stackrel{\mathrm
{a.s.}}{\rightarrow}F(u)$ uniformly in $u$, where $F(u)$ is the
C.D.F. of the pooled sample.

\item[(A4)] $\eta_{\min}\equiv\min_{1\leq r\leq
K_n}\int_0^1\eta_r(u)/\{F(u)(1-F(u))\}\,dF(u)$ is a positive constant.
\end{longlist}
Assumption (A1) is required in some exponential tail inequalities as
detailed in the proof of Lemma \ref{klem4}, while the $F_k$'s can be
discrete or mixed distributions in practice. Assumption (A2) is
a standard requirement for the theoretical development in the MCP,
which allows the
change-points to be asymptotically distinguishable. Assumption (A3)
is a technical condition that is trivially satisfied
by the Glivenko--Cantelli theorem when $K_n$ is finite.\vspace*{1pt} Generally, it
can be
replaced by the conditions that\vspace*{1pt}
$\lim_{n\rightarrow\infty}\sum_{k=1}^{K_n+1}(\tau_k-\tau_{k-1})/nF_k(u)$
exists and
$\sum_{k=1}^{K_n+1} \{(\tau_k-\tau_{k-1})/n\sup_{u}|{\widehat
{F}}_{\tau
_{k-1}}^{\tau_k}(u)-F_k(u)| \}$
converges to 0 a.s. By the Dvoretzky--Kiefer--Wolfowitz inequality,
the latter one holds if
$\sum_{n=1}^{\infty}K_n\exp(-2\lambda_n\epsilon^2/K_n^2)<\infty$ for
any $\epsilon>0$. Assumption (A4) means that the smallest signal
strength among all the changes is bounded away from zero.

We may consider relaxing $\lambda_n\rightarrow\infty$ in assumption~(A2)
by allowing $\eta_{\min}\rightarrow\infty$ as $n\rightarrow\infty
$. It is
intuitive that if two successive distributions are very different, then
we do not need a very large $\lambda_n$ to locate the change point.
For the mean change problem, \citet{NiuZha12} and \citet{HaoNiuZha13} revealed that in order to obtain the $O_p(1)$ consistency, a
condition $\delta\lambda_n>32\log n$ is required, where $\delta$ is
the minimal jump size at the change-points (similar to $\eta_{\min}$).
In our nonparametric setting, such an extension warrants future investigation.

%
\begin{teo}\label{teo1}
Under assumptions \textup{(A1)--(A4)}, if $K_n^3(\log
K_n)^2(\log\delta_n)^2/\delta_n$ $\rightarrow0$ and
$\delta_n/\lambda_n\rightarrow0$, then
\[
\Pr\bigl\{\mathcal{G}_n(K_n)\in C_{K_n}(
\delta_n) \bigr\} \rightarrow1\qquad\mbox{as } n\rightarrow\infty.
\]
\end{teo}

Under the classical mean change-point model, \citet{YaoAu89}
studied the property of the least-squares estimator,
%
%
\begin{equation}
\label{yao} \mathop{\arg\min}_{\tau_1'<\cdots<\tau_{K_n}'}\sum_{k=1}^{K_n+1}
\sum_{i=\tau_{k-1}'}^{\tau_k'-1}\bigl\{X_i-
\hat{\mu}\bigl(\tau_{k-1}',\tau_k'
\bigr)\bigr\}^2,
\end{equation}
where $\hat{\mu}(\tau_{k-1}',\tau_k')$ denotes the average of the
observations $\{X_{\tau_{k-1}'},\ldots,X_{{\tau_k'}-1}\}$. It is
well known that the least-squares estimator is consistent
with the optimal rate $O_p(1)$, when the number of change-points is
known (and does not depend on $n$) and the change magnitudes are
fixed; see \citet{HaoNiuZha13} and the references therein. Under a
similar setting with $K_n\equiv K$, we can establish the same rate
of $O_p(1)$ for our nonparametric approach.

%
\begin{coro}\label{cor1}
Under assumptions \textup{(A1)}, \textup{(A2)} and \textup{(A4)},
$|\hat{\tau}_s-\tau_s|=O_p(1)$ for $s=1,\ldots, K$.
\end{coro}
The proof is similar to that of Theorem~\ref{teo1},
which is provided in the supplementary material [\citet{Zouetal}].
With the knowledge of $K$, we can obtain an optimal rate of $O_p(1)$
without specifying the distributions, which is
consistent with the single change-point case in \citet{Dum91}.

The next theorem establishes the consistency of the NMCD procedure
with the BIC in (\ref{BIC}).
Let $\widehat{K}_n=\arg\min_{1\leq L\leq\widebar{K}_{n}}{\mathrm
{BIC}}_L$, where
$\widebar{K}_{n}$ is an upper bound on the true number of change-points.
%
%
\begin{teo}\label{probic}
Under\vspace*{1pt} assumptions \textup{(A1)--(A4)}, $\lambda_n/(\widebar{K}_{n}\zeta_n)
\to
\infty$, $\zeta_n= \widebar{K}_{n}^3(\log\widebar{K}_{n})^2
(\log
n )^{2+c}$ with any $c>0$, then $\Pr(\widehat
{K}_n=K_n)\rightarrow
1$ as
$n\rightarrow\infty$.
\end{teo}
It is remarkable that in the conventional setting where
$\widebar{K}_{n}$ is bounded, we can use $\zeta_n$ of order $(\log
n)^{2+c}$ instead of its least-squares counterpart $\log n$ in
\citet{Yao88}. In conjunction with Theorem~\ref{teo1}, this result implies that
$\Pr\{\mathcal{G}_n(\widehat{K})\in C_{K}((\log n)^{2+c})\}
\rightarrow
1$ with a fixed number of change-points.

\section{Implementation of NMCD}\label{sec3}

\subsection{Algorithm}\label{sec3.1}

One important property of the proposed maximum likelihood
approach is that (\ref{fpl}) is separable. The optimum for splitting
cases $1,\ldots,n$ into~$L$ segments conceptually consists of first
finding the rightmost change-point $\hat{\tau}_L$, and then finding
the remaining change-points from the fact that they constitute the
optimum for splitting cases $1,\ldots,\hat{\tau}_L$ into $L-1$
segments. This separability is called Bellman's ``principle of
optimality'' [\citet{BelDre62}]. Thus, (\ref{fpl}) can be
maximized via the DP algorithm and fitting such a nonparametric MCP
model is straightforward and fast. The total computational
complexity is $O(Ln^2)$ for a given $L$; see \citet{Haw01} and
\citet{BaiPer03} for the pseudo-codes of the DP. \citet{Haw01}
suggested using the DP on a grid of $m\ll n$ values. \citet
{HarLev10} proposed using a LASSO-type
penalized estimator to achieve a reduced version of the
least-squares method. \citet{NiuZha12} developed a screening and
ranking algorithm to detect DNA copy number variations in the MCP
framework.

Due to the DP's computational complexity in $n^2$,
an optimal segmentation of a
very long sequence could be computationally intensive; for example,
DNA sequences nowadays are often extremely long
[\citet{FeaVas09}].
To alleviate the computational burden,
we introduce a preliminary screening step which can
exclude most of the irrelevant points and, as a consequence,
the NMCD is implemented in a much lower-dimensional space.

\subsubsection*{Screening algorithm}

\begin{longlist}[(iii)]
\item[(i)] Choose an appropriate integer $n_{I}$ which is the length
of each subsequence of the data, and take the
estimated change-point set $\mathcal{O}=\varnothing$.
\item[(ii)]
Initialize $\gamma_i=0$ for $i=1,\ldots,n$; and for $i=n_I,\ldots,n-n_I$,
update $\gamma_i$ to be the Cram\'{e}r--von Mises two-sample test
statistic for the samples
$\{X_{i-n_I+1},\ldots,X_i\}$ and $\{X_{i+1},\ldots, X_{i+n_I}\}$.

\item[(iii)] For $i=n_I,\ldots,n-n_I$, define $k=\mathop{\arg\max
}_{i-n_I< j\leq
i+n_I}\gamma_j$. If $k=i$, update $\mathcal{O}=\mathcal{O}\cup\{i\}$.
\end{longlist}

Intuitively speaking, this screening step finds the most influential
points that have the largest \textit{local} jump sizes quantified by
the Cram\'{e}r--von Mises statistic, and thus helps to avoid
including too many candidate points around the true change-point. As
a result, we can obtain a candidate change-point set,
$\mathcal{O}$, of which the cardinality, $|\mathcal{O}|$, is usually
much smaller than $n$. Finally, we run the NMCD procedure within the\vadjust{\goodbreak}
set $\mathcal{O}$ using the DP algorithm to find the solution of
\[
\mathop{\arg\max}_{\tau_1'<\cdots<\tau_L'\in
\mathcal{O}}R_n\bigl(\tau_1',
\ldots,\tau_L'\bigr).
\]

Apparently, the screening procedure is fast because it
mainly requires calculating $n-2n_I+1$ Cram\'{e}r--von Mises statistics.
In contrast, \citet{Lee96} used a thresholding step to determine the
number of change-points. The main difference between \citet{Lee96}
and \citet{NiuZha12} lies in the choice
of the local test statistic; the former uses some seminorm of
empirical distribution functions and the latter is based on the
two-sample mean difference.

We next clarify how to choose $n_I$, which
formally establishes the consistency of the screening procedure.

%
\begin{pro}\label{SIS} Under assumptions \textup{(A1)--(A2)},
if $n_I/\log n\rightarrow\infty$ and $n_I/\lambda_n^{1/2}\rightarrow
0$, then we
have $\Pr\{\mathcal{O}\in H_{|\mathcal{O}|}(\log n) \}
\rightarrow
1$, where
\begin{eqnarray*}
H_{l}(\delta_n)&=& \bigl\{\bigl(\tau_1',
\ldots,\tau_{l}'\bigr)\dvtx  1<\tau_1'<
\cdots<\tau_{l}'\leq n, \mbox{ and for each } 1\leq r\leq
K_n
\\
&&\hspace*{73pt} \mbox{there exists at least a } \tau_s'
\mbox{ so that } |\tau_s'-\tau_r|\leq
\delta_n \bigr\}.
\end{eqnarray*}
\end{pro}
This result follows by
verifying condition (A3) in \citet{Lee96}; see Example~II
of \citet{Dum91}.
With probability tending to one, the screening algorithm
can at least include one $\delta_n$-neighborhood of the true
location set by choosing an appropriate $n_I$. Given a candidate
$L$, the computation of NMCD reduces to $O(L|\mathcal{O}|n)$,
which is of order $O(\widebar{K}_{n}^2|\mathcal{O}|n)$ in conjunction
with the BIC. Both the R and \mbox{FORTRAN} codes for
implementing the entire procedure are available from the authors upon
request.

\subsection{Selection of tuning parameters}\label{sec3.2}
We propose to take $dw(u)=\break \{{\widehat{F}}_n(u)(1-{\widehat{F}}_n(u))\}
^{-1}\,d{\widehat{F}}_n(u)$, which\vspace*{1pt}
is found to be more powerful than simply using
$dw(u)=d{\widehat{F}}_n(u)$. The function
$\{{\widehat{F}}_n(u)(1-{\widehat{F}}_n(u))\}^{-1}$ attains its
minimum at ${\widehat{F}}_n(u)=1/2$,
that is when $u$ is the median of the sample. Intuitively, when two
successive distributions mainly differ in their centers, both
choices of $dw(u)$ would be powerful \mbox{because} a large portion of
observations are around the center. However, if the difference
between two adjacent distributions lies in their tails, using
$dw(u)=d{\widehat{F}}_n(u)$ may not work well because only very
limited information
is included in the integral of (\ref{pl}). In contrast,
our weight would be larger for those more extreme observations
(far way from the median).

To better understand this, we analyze the term $\eta_{\min}$,
which reflects the detection ability to a large extent.
Consider a special case
\begin{eqnarray*}
X_i \sim\cases{ U(0,1), &\quad$1\leq i \leq{n}/{2}$,
\vspace*{3pt}\cr
U(1,2), &
\quad${n}/{2}+1\leq i \leq n$}\vadjust{\goodbreak}
\end{eqnarray*}
and thus
\begin{eqnarray*}
\eta_1(u)&=& \biggl(u \log2+(1-u)\log\frac{1-u}{1-{u}/{2}} \biggr)
I(0<u<1)
\\
&&{}+ \biggl((u-1)\log\frac{2(u-1)}{u}+(2-u)\log\frac{2-u}{1-{u}/{2}}
\biggr)I(1<u<2).
\end{eqnarray*}
It is easy to check that $\eta_{\min}=\int_0^2
{\eta_1(u)}/\{F(u)(1-F(u))\}\,dF(u)$ is unbounded, while the counterpart
$\int_0^2 \eta_1(u)\,dF(u)$ is finite.
Consequently,\vspace*{1pt} the NMCD\break  procedure would be more powerful
by using the weight $\{{\widehat{F}}_n(u)(1-\break {\widehat{F}}_n(u))\}
^{-1}\,d{\widehat{F}}_n(u)$.\vspace*{2pt}

Under the assumption that $\zeta_n=\widebar{K}_{n}^3(\log
\widebar{K}_{n})^2 (\log n )^{2+c}$ with $c>0$ and
$\lambda_n/(\widebar{K}_{n}\zeta_n) \to\infty$, we establish the
consistency of the BIC in (\ref{BIC}) for model selection.
The choice of $\zeta_n$ depends on
$\widebar{K}_{n}$ and $\lambda_n$ which are unknown. The value of
$\widebar{K}_{n}$ depends on the practical consideration of how many
change-points are to be identified, while $\lambda_n$ reflects the
length of the smallest segment.
For practical use, we take $\widebar{K}_{n}$ to be fixed
and recommend $\zeta_n=(\log n)^{2+c}/2$ with $c=0.1$.
A~small value of $c$ helps to prevent
underfitting, as one is often reluctant to miss any important
change-point. The\vspace*{1pt} performance of NMCD insensitive to
the choice of $\widebar{K}_{n}$, as long as $\widebar{K}_{n}$ is not
too small,
which is also to avoid underfitting. We suggest
$\widebar{K}_{n}=|\mathcal{O}|$, that is, the cardinality of the candidate
change-point set in the screening algorithm.

\section{Simulation studies}\label{sec4}
\subsection{Model setups}\label{sec4.1}
To evaluate the finite-sample performance of the proposed NMCD
procedure, we
conduct extensive simulation studies, and also make comparisons with
existing methods. We calculate the distance
between the estimated set $\widehat{\mathcal{G}}_n$ and the true
change-point set $\mathcal{C}_t$ [\citet{Boyetal09}],
\[
\xi(\widehat{\mathcal{G}}_n\|\mathcal{C}_t)=\sup
_{b\in
\mathcal{C}_t} \inf_{a\in\widehat{\mathcal{G}}_n}|a-b|\quad\mbox{and}\quad
\xi(\mathcal{C}_t\|\widehat{\mathcal{G}}_n)=\sup
_{b\in\widehat
{\mathcal{G}}_n} \inf_{a\in\mathcal{C}_t}|a-b|,
\]
which quantify the over-segmentation error and
the under-segmentation error, respectively. A desirable estimator should
be able to balance both quantities. In addition, we consider the
average Rand index [\citet{FowMal83}], which
measures the discrepancy of two sets from an average viewpoint.

Following model~(I) introduced by \citet{DonJoh95}, we
generate the Blocks datasets, which contains $K_n=11$ change-points:
\begin{eqnarray*}
\mbox{Model (I)}\dvt X_i&=&\sum_{j=1}^{K_n}
h_jJ(nt_i-\tau_j)+\sigma{
\varepsilon}_i,\qquad J(x)=\bigl\{1+\operatorname{sgn}(x)\bigr\}/2,
\\
\{\tau_j/n\}&=&\{0.1,0.13,0.15,0.23,0.25,0.40,0.44,0.65,0.76,0.78,0.81
\},
\\
\{h_j\}&=&\{2.01,-2.51,1.51,-2.01,2.51,-2.11,1.05,2.16,
\\
&&\hspace*{135pt}{} -1.56,2.56,-2.11\},
\end{eqnarray*}
where there are $n$ equally spaced covariates $t_i$ in $[0,1]$.
Three error distributions for~${\varepsilon}_i$ are considered: $N(0,1)$,
Student's $t$ distribution with three degrees of freedom $t_{(3)}$, and
the standardized (zero mean and unit variance) chi-squared
distribution with one degree of freedom $\chi^2_{(1)}$. The Blocks
datasets with $n=1000$, as depicted in the top three plots of Figure~A.1 in the supplementary material [\citet{Zouetal}], are generally considered difficult for
multiple change-point estimation due to highly heterogeneous segment
levels and lengths.

In a more complicated setting with both location and scale changes,
we consider model~(II) with $K_n=4$:
\begin{eqnarray*}
\mbox{Model (II)}\dvt X_i&=&\sum_{j=1}^{K_n}
h_jJ(nt_i-\tau_j)+\sigma{
\varepsilon}_i \prod_{j=1}^{\sum_{j=1}^{K_n}{J(nt_i-\tau_j)}}
v_j,
\\
\{h_j\}&=&\{3,0,-2,0\},\qquad
\{\tau_j/n\}=\{0.20,0.40,0.65,0.85\}\quad\mbox{and}
\\
\{v_j\}&=&\{1,5,1,0.25\},
\end{eqnarray*}
where all the other setups are the same as those of model~(I). As shown
by the bottom three plots in Figure A.1, there are two location
changes and two scale changes.

In addition, we include a simulation study when the distributions
differ in the skewness and kurtosis. In particular, we consider
\[
\mbox{Model (III)}\dvt X_i\sim F_j(x), \qquad
\tau_j/n=\{0.20,0.50,0.75\}, \qquad j=1,2,3,4,
\]
where $F_1(x), \ldots, F_4(x)$ correspond to the standard normal,
the standardized $\chi^2_{(3)}$ (with zero mean and unit variance), the
standardized $\chi^2_{(1)}$, and the standard normal distribution,
respectively. Because there is no mean or variance difference
between the $F_j$'s, as depicted in the left panel of Figure A.4, the
estimation for such a change-point problem is rather difficult. All
the simulation results are obtained with 1000 replications.

\subsection{Calibration of tuning parameters}\label{sec4.2}
\begin{figure}

\includegraphics{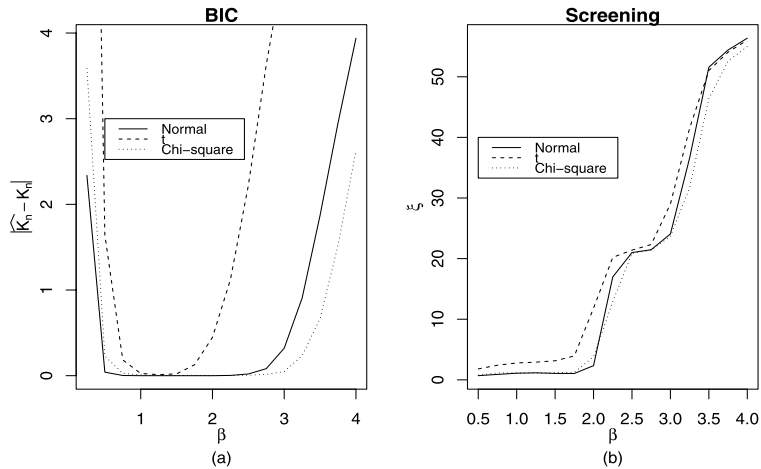}

\caption{The performance of NMCD under model~\textup{(I)} with $n=1000$ and $\sigma=0.5$ when the tuning parameters vary:
\textup{(a)}~the curves of $|\widehat{K}_n-K_n|$ versus the value of $\beta$;
\textup{(b)} the curves of $\xi(\widehat{\mathcal{G}}_n\|\mathcal{C}_t)$ versus
the value of $\beta$.}\label{fig1}
\end{figure}

To study the sensitivity of the choice of $\zeta_n$, Figure~\ref{fig1}(a)
shows the
curves of $|\widehat{K}_n-K_n|$ versus the value of $\beta$ with
$\zeta_n=\beta(\log n)^{2.1}/2$ under model~(I).
Clearly, the estimation is reasonably well with a value of $\beta$
around 1.
For more adaptive model selection,
a data-adaptive complexity penalty in \citet{SheYe02} could be considered.

%

In the screening procedure, the choice of $n_I$ needs to balance the
computation and underfitting.
By Proposition \ref{SIS}, $n_I\in(\log n,\lambda_n^{1/2})$,
while\vadjust{\goodbreak}
$\lambda_n$ is typically unknown.
In practice, we recommend to choose $n_I=\lceil(\log
n)^{3/2}/2\rceil$, which is the smallest integer that is larger than
$(\log n)^{3/2}/2$. Figure~\ref{fig1}(b) shows the curves of
under-segmentation errors versus the value of $\beta$ with
$n_I=\lceil\beta(\log n)^{3/2}/2\rceil$ under model~(I). In a
neighborhood of $\beta=1$, our method provides a reasonably
effective reduction of the subset $\mathcal{O}$ and the performance
is relatively stable. In general, we do not recommend
a too large value of $n_I$ so as to avoid underfitting.
From the results shown in Section~\ref{sec4.6}, the choice of $\zeta
_n=(\log
n)^{2.1}/2$ and $n_I=\lceil(\log n)^{3/2}/2\rceil$ works also well
when the number of change-points increases as the sample size increases.

%
\begin{table}
\tabcolsep=0pt
\caption{Comparison of the parametric likelihood (PL), NMCD, and NMCD* methods when the number
of change-points $K_n$ is specified (known) under models \textup{(I)}~and~\textup{(II)}, respectively.
The standard deviations are given in parentheses}\label{TknowK}
\begin{tabular*}{\tablewidth}{@{\extracolsep{\fill}}@{}lcd{4.0}cccccc@{}}
\hline
& & & \multicolumn{3}{c}{$\bolds{\xi(\widehat{\mathcal{G}}_n\|\mathcal{C}_t)}$} & \multicolumn{3}{c@{}}{$\bolds{\xi(\mathcal{C}_t\|\widehat{\mathcal{G}}_n)}$}\\[-6pt]
& & & \multicolumn{3}{c}{\hrulefill} & \multicolumn{3}{c@{}}{\hrulefill}\\
\textbf{Model} & \textbf{Error} & \multicolumn{1}{c}{$\bolds{n}$}&
\multicolumn{1}{c}{\textbf{PL}} & \multicolumn{1}{c}{\textbf{NMCD}} & \multicolumn{1}{c}{\textbf{NMCD*}}
& \multicolumn{1}{c}{\textbf{PL}} & \multicolumn{1}{c}{\textbf{NMCD}} & \textbf{NMCD*} \\
\hline
(I)&$N (0,1)$& 500 & 0.96\ (1.19) & 0.96\ (1.14) & 1.16\ (1.15) & 0.96\ (1.19)& 0.96\ (1.14) & 1.16\ (1.15) \\
& & 1000 & 0.91\ (1.15) & 0.97\ (1.16) & 1.06\ (1.21) & 0.91\ (1.15) &0.97\ (1.16) & 1.06\ (1.21) \\
&$t_{(3)}$& 500 & 13.6\ (12.0) & 3.77\ (4.48) & 3.86\ (4.33) & 14.3\ (18.4) &3.95\ (7.51) & 3.97\ (7.63) \\
& & 1000 & 20.2\ (21.3) & 2.58\ (2.50) & 2.90\ (2.72) & 21.9\ (34.5) &2.56\ (2.40) & 2.90\ (2.72) \\
&$\chi^2_{(1)}$ & 500 & 1.39\ (2.91) & 0.70\ (0.80) & 0.80\ (1.22) &1.13\ (1.57) & 0.70\ (0.80) & 0.81\ (1.41) \\
& & 1000 & 1.05\ (2.15) & 0.59\ (0.77) & 0.58\ (0.71) & 0.99\ (1.38) &0.59\ (0.77) & 0.58\ (0.71)
\\[3pt]
(II)&$N (0,1)$&500 & 1.59\ (1.72) & 2.35\ (2.42) & 3.34\ (4.96) & 1.59\ (1.72)& 2.35\ (2.42) & 3.34\ (4.96) \\
&& 1000 & 1.58\ (1.52) & 2.68\ (2.59) & 2.74\ (2.89) & 1.58\ (1.52) &2.68\ (2.59) & 2.74\ (2.89) \\
&$t_{(3)}$ &500 & 13.6\ (25.8) & 4.75\ (6.87) & 6.42\ (8.84) & 7.52\ (10.2) &4.54\ (5.19) & 6.05\ (6.42) \\
&& 1000 &16.4\ (40.2) & 4.10\ (3.88) & 5.27\ (7.20) & 10.3\ (18.0) &4.10\ (3.88) & 5.24\ (6.85) \\
&$\chi^2_{(1)}$ &500 &6.36\ (11.3) & 1.57\ (2.12) & 1.65\ (2.90) &5.88\ (8.93) & 1.57\ (2.12) & 1.65\ (2.90) \\
&& 1000 &4.80\ (67.8) & 1.17\ (1.45) & 1.49\ (2.10) & 4.80\ (7.82) &1.17\ (1.45) & 1.49\ (2.10) \\
\hline
\end{tabular*}
\end{table}

\subsection{Comparison between NMCD and PL}\label{sec4.3}
Firstly, under model~(I) with location changes only, we make a
comparison of NMCD with the parametric likelihood (PL) method which
coincides with the classical least-squares method in (\ref{yao})
under the normality assumption [\citet{Yao88}]. We also consider a
variant of NMCD by using $dw(u)=d{\widehat{F}}_n(u)$ (abbreviated as NMCD*).
The comparison is conducted with and without knowing the true number
of change-points $K_n$, respectively. Table~\ref{TknowK} presents
the average values of $\xi(\widehat{\mathcal{G}}_n\|\mathcal{C}_t)$ and
$\xi(\mathcal{C}_t\|\widehat{\mathcal{G}}_n)$ for $n=500$ and 1000 and
$\sigma=0.5$ when $K_n$ is known to be 11. To gain more insight,
we also present the standard deviations of the two distances in
parentheses. Simulation results with other values of
$\sigma$ can be found in the supplementary material [\citet{Zouetal}].

As expected, the PL has superior efficiency for the case with
normal errors, since the parametric model is correctly specified.
The NMCD procedure also offers satisfactory performance and the
differences in the two $\xi$ values between NMCD and PL are
extremely small, while both methods significantly outperform the
NMCD* procedure. For the cases with $t_{(3)}$ and $\chi^2_{(1)}$ errors,
the NMCD procedure almost uniformly outperforms the PL in terms of
estimation accuracy of the locations. Not only are the distance
values of $\xi(\widehat{\mathcal{G}}_n\|\mathcal{C}_t)$ and
$\xi(\mathcal{C}_t\|\widehat{\mathcal{G}}_n)$ smaller, but the
corresponding standard deviations are also much smaller using
the NMCD.

%
\begin{table}
\tabcolsep=0pt
\caption{Comparison of the PL and NMCD methods
when the number of change-points $K_n$ is unknown ($K_n$~is selected
using the BIC) under models \textup{(I)} and \textup{(II)}, respectively.
The standard deviations are given in parentheses}\label{TunknowK}
\begin{tabular*}{\tablewidth}{@{\extracolsep{\fill}}@{}lcd{4.0}cccccc@{}}\hline
& & & \multicolumn{3}{c}{\textbf{Parametric likelihood (PL)}} & \multicolumn{3}{c@{}}{\textbf{NMCD}}\\[-6pt]
& & & \multicolumn{3}{c}{\hrulefill} & \multicolumn{3}{c@{}}{\hrulefill}
\\
\textbf{Model} & \textbf{Error} & \multicolumn{1}{c}{$\bolds{n}$} & \multicolumn{1}{c}{$\bolds{\xi(\widehat{\mathcal{G}}_n\|\mathcal{C}_t)}$} &
\multicolumn{1}{c}{$\bolds{\xi(\mathcal{C}_t\|\widehat{\mathcal{G}}_n)}$} & \multicolumn{1}{c}{$\bolds{|\widehat{K}_n-K_n|}$} &
\multicolumn{1}{c}{$\bolds{\xi(\widehat{\mathcal{G}}_n\|\mathcal{C}_t)}$} & \multicolumn{1}{c}{$\bolds{\xi(\mathcal{C}_t\|\widehat{\mathcal{G}}_n)}$} &
$\bolds{|\widehat{K}_n-K_n|}$\\
\hline
(I)&$N (0,1)$& 500 & 0.93\ (1.08) & 2.16\ (6.57) & 0.09\ (0.31) & 0.96\ (1.34)& 0.99\ (1.05) & 0.00\ (0.04) \\
& & 1000 & 0.94\ (1.14) & 2.30\ (10.3) & 0.05\ (0.25) & 0.96\ (1.25) &1.01\ (1.25) & 0.00\ (0.04) \\
&$t_{(3)}$& 500 & 2.91\ (2.92) & 39.0\ (24.9) & 6.05\ (3.47) & 3.34\ (4.22) &8.64\ (15.2) & 0.36\ (0.88) \\
& & 1000 & 2.94\ (3.02) & 95.2\ (48.8) & 9.70\ (4.14) & 2.54\ (2.78) &10.0\ (26.8) & 0.36\ (0.75) \\
&$\chi^2_{(1)}$ & 500 & 0.85\ (0.99) & 49.5\ (23.6) & 10.9\ (4.69) &0.73\ (0.95) & 1.36\ (5.59) & 0.05\ (0.28) \\
& & 1000 & 0.85\ (1.05) & 111\ (46.2)  & 14.2\ (4.06) & 0.53\ (0.69) &0.89\ (4.28) & 0.02\ (0.20)
\\[3pt]
(II)& $N (0,1)$ & 500 & 1.66\ (1.61) & 2.22\ (5.56) & 0.04\ (0.22) &2.28\ (2.31) & 4.45\ (8.54) & 0.13\ (0.37) \\
& & 1000 & 1.69\ (1.50) & 1.71\ (1.52) & 0.01\ (0.11) & 2.19\ (2.11) &3.93\ (10.6) & 0.06\ (0.27) \\
& $t_{(3)}$ & 500 & 5.77\ (6.57) & 24.1\ (20.0) & 1.58\ (1.56) & 5.18\ (6.18)& 14.1\ (16.5) & 0.75\ (1.01)\\
& & 1000 & 5.59\ (6.26) & 62.4\ (41.3) & 2.72\ (2.21) & 4.50\ (4.44) &17.0\ (28.4) & 0.47\ (0.87) \\
& $\chi^2_{(1)}$ & 500 & 5.03\ (6.19) & 43.1\ (16.0) & 4.71\ (2.66) &1.67\ (2.39) & 7.27\ (12.6) & 0.43\ (0.80) \\
& & 1000 & 5.00\ (6.29) & 91.1\ (31.1) & 6.22\ (3.23) & 1.26\ (1.50) &9.45\ (22.7) & 0.28\ (0.70) \\
\hline
\end{tabular*}
\end{table}

Next, we consider the $K_n$ unknown case, for which both the NMCD
and PL procedures are implemented by setting $\widebar{K}_{n}=30$ and using
the BIC to choose the number of change-points. The average values of
the distances $\xi(\widehat{\mathcal{G}}_n\|\mathcal{C}_t)$ and
$\xi(\mathcal{C}_t\|\widehat{\mathcal{G}}_n)$ are tabulated in
Table~\ref{TunknowK}. In addition,\vspace*{2pt} we also present the average values of
$|\widehat{K}_n-K_n|$ with standard deviations in parentheses, which
reflect the overall estimation accuracy of $K_n$. Clearly, the two
methods have comparable performances under the normal error, while
the proposed NMCD significantly outperforms PL in terms of
$\xi(\mathcal{C}_t\|\widehat{\mathcal{G}}_n)$ and $|\widehat
{K}_n-K_n|$ for
the two nonnormal cases, because the efficiency of the BIC used in
PL relies heavily on the parametric assumption. When we compare the
results across Tables~\ref{TknowK} and \ref{TunknowK}, the standard deviations for the
distance measures increase from the $K_n$ known to the $K_n$ unknown
cases, as estimating $K_n$ further enlarges the variability.

We turn to the comparison between NMCD and PL under model (II) in
which both location and scale changes are exhibited. In this
situation, the standard least-squares method (\ref{yao}) does not
work well because it is constructed for location changes only. To
further allow for scale changes under the PL method, we consider
%
%
\begin{equation}
\label{yaoprime} \mathop{\arg\min}_{\tau_1'<\cdots<\tau_{K_n'}}\sum
_{k=1}^L \bigl(\tau_{k+1}'-
\tau_k'\bigr)\log\hat{\sigma}_k^2,
\end{equation}
where $\hat{\sigma}_k^2=(\tau_{k+1}'-\tau_k')^{-1}
\sum_{i=\tau_{k-1}'}^{\tau_k'-1}\{X_i-\hat{\mu}(\tau_{k-1}',\tau
_k')\}^2$,
and the BIC is modified accordingly. The\vspace*{1pt} bottom panels of Tables~\ref
{TknowK} and \ref{TunknowK} tabulate the values of
$\xi(\widehat{\mathcal{G}}_n\|\mathcal{C}_t)$\vadjust{\goodbreak} and
$\xi(\mathcal{C}_t\|\widehat{\mathcal{G}}_n)$ when $K_n$ is
specified in
advance and estimated by using the BIC, respectively. Clearly, the
NMCD method delivers a satisfactory detection performance for the
normal case and performs much better than the PL method for the two
nonnormal cases. Therefore, the conclusion remains that the PL
method is generally sensitive to model specification, while the NMCD
does not depend on any parametric modeling assumption and thus is
much more robust.

\subsection{Comparisons of NMCD with other nonparametric methods}\label{sec4.4}
We consider the methods of \citet{Lee96} and \citet{MatJam},
as they also do not make any assumptions regarding the nature of the changes.
The NMCD is implemented with the initial nonparametric screening
procedure, and $K_n$ is selected by the BIC. In both our screening
procedure and Lee's (\citeyear{Lee96}) method, the window is set as
$n_{I}=\lceil(\log n)^{3/2}/2\rceil$, and the threshold value of the
latter is chosen as $(\log n)^{3/4}$. The ECP method of
\citet{MatJam} is implemented
using the ``ecp'' R package with the false alarm rate 0.05 and
$\alpha=1$.\vspace*{1pt}\vadjust{\goodbreak}

%
\begin{table}
\tabcolsep=0pt
\caption{Comparison of NMCD, Lee's (\citeyear{Lee96}) method and
Matteson and James's (\citeyear{MatJam}) ECP in terms of $\xi(\widehat{\mathcal{G}}_n;\mathcal
{C}_t)$, Rand and $|\widehat{K}_n-K_n|$ under models \textup{(I)--(III)} with $\sigma=0.5$}\label{TLeeunknowK}
\begin{tabular*}{\tablewidth}{@{\extracolsep{\fill}}@{}lcd{4.0}d{3.1}d{3.2}d{2.2}cccccc@{}}
\hline
& & &\multicolumn{3}{c}{$\bolds{\xi(\widehat{\mathcal{G}}_n;\mathcal{C}_t)}$} & \multicolumn{3}{c}{\textbf{Rand}} & \multicolumn{3}{c@{}}{$\bolds{|\widehat{K}_n-K_n|}$}\\[-6pt]
& & &\multicolumn{3}{c}{\hrulefill} & \multicolumn{3}{c}{\hrulefill} & \multicolumn{3}{c@{}}{\hrulefill}
\\
\textbf{Model} & \textbf{Error} & \multicolumn{1}{c}{$\bolds{n}$} & \multicolumn{1}{c}{\textbf{Lee}}
& \multicolumn{1}{c}{\textbf{ECP}} & \multicolumn{1}{c}{\textbf{NMCD}} & \multicolumn{1}{c}{\textbf{Lee}}
& \multicolumn{1}{c}{\textbf{ECP}} & \multicolumn{1}{c}{\textbf{NMCD}} & \multicolumn{1}{c}{\textbf{Lee}}
& \multicolumn{1}{c}{\textbf{ECP}} & \textbf{NMCD}
\\
\hline
(I) & $N(0,1)$ & 500 & 84.9 & 6.03 & 2.62 & 0.920 & 0.994 & 0.992 &28.5 & 0.07 & 0.01 \\
& & 1000 & 176 & 7.42 & 2.23 & 0.915 & 0.997 & 0.994 & 43.2 & 0.07 &0.00 \\
& $t_{(3)}$ & 500 & 86.7 & 4.95 & 8.94 & 0.920 & 0.995 & 0.988 & 27.5& 0.06 & 0.22 \\
& & 1000 & 177 & 7.29 & 7.63 & 0.914 & 0.997 & 0.993 & 42.9 & 0.08 &0.02 \\
& $\chi^2_{(1)}$ & 500 & 85.0 & 4.67 & 3.00 & 0.921 & 0.995 & 0.992& 28.3 & 0.06 & 0.02 \\
& & 1000 & 176 & 5.67 & 2.80 & 0.915 & 0.997 & 0.994 & 43.1 & 0.05 & 0.01
\\[3pt]
(II) & $N(0,1)$ & 500 & 69.0 & 17.6 & 14.4 & 0.832 & 0.980 & 0.980 &33.8 & 0.06 & 0.11 \\
& & 1000 & 140 & 17.5 & 14.4 & 0.830 & 0.990 & 0.987 & 51.3 & 0.07 &0.03 \\
& $t_{(3)}$ & 500 & 69.3 & 16.8 & 20.4 & 0.833 & 0.982 & 0.974 & 33.7& 0.10 & 0.25 \\
& & 1000 & 141 & 12.5 & 21.4 & 0.830 & 0.992 & 0.983 & 51.6 & 0.06 &0.13 \\
& $\chi^2_{(1)}$ & 500 & 67.9 & 8.25 & 10.5 & 0.833 & 0.989 & 0.983& 34.0 & 0.05 & 0.12 \\
& & 1000 & 139 & 10.2 & 12.6 & 0.830 & 0.994 & 0.987 & 51.2 & 0.07 &0.09
\\[3pt]
(III) & & 500 & 120 & 394 & 78.2 & 0.822 & 0.446 & 0.894 & 35.4 &1.73 & 0.53 \\
& & 1000 & 243 & 452 & 43.9 & 0.818 & 0.714 & 0.965 & 52.8 & 1.22 &0.19 \\
\hline
\end{tabular*}
\end{table}

Table~\ref{TLeeunknowK} shows the comparison results based on
$\xi(\widehat{\mathcal{G}}_n;\mathcal{C}_t)\equiv\xi(\widehat
{\mathcal{G}}_n\|\mathcal{C}_t)+
\xi(\mathcal{C}_t\|\widehat{\mathcal{G}}_n)$, $|\widehat
{K}_n-K_n|$, and
the Rand index
under models (I)--(III) with $\sigma=0.5$, respectively.
Lee's (\citeyear{Lee96}) method is unable to produce a reasonable
estimate for
$K_n$ and the
resulting models are much overfitted in all the cases, which indicates
that its ``local'' nature incurs substantial loss of the information.
Under model~(I), the NMCD performs better than ECP for normal\vspace*{2pt} and $\chi
^2_{(1)}$ errors,
while the opposite is true for the $t$ error distribution.
Under model (II), the ECP also exhibits certain advantage, especially
for Student's $t$
and $\chi^2_{(1)}$ error distributions.
Both the NMCD and ECP methods significantly outperform that of
\citet{Lee96} in
models (I) and (II). Under model (III), both the methods of
ECP and \citet{Lee96} appear not working well, while
the NMCD still produces reasonable detection results.
As the divergence measure used in the ECP is
essentially similar to Euclidean distances, the ECP is
expected to perform well when the distributions differ in
the first two moments, which however is not the case for model (III).
The advantages of NMCD are mainly due to the joint use of the
nonparametric likelihood and the weight function
$w(u)=\{{\widehat{F}}_n(u)(1-{\widehat{F}}_n(u))\}^{-1}\,d{\widehat{F}}_n(u)$.
Based on the empirical distribution functions,
the nonparametric likelihood approach is capable of detecting
various types of changes. In addition, the difference between two
adjacent distributions under model (III) does not
lie in their centers, and thus using our proposed $w(u)$ would provide certain
improvement as discussed in Section~\ref{sec3.2}.
Due to the use of DP, our procedure is much faster than the ECP.

\subsection{Comparison of NMCD and LSTV}\label{sec4.5}
\citet{HarLev10} proposed
the least-squares total variation method (LSTV) to estimate the
locations of multiple change-points. By reframing the MCP
in a variable selection context, they use a penalized least-squares
criterion with a LASSO-type penalty. The LSTV enjoys efficient
computation using the least angle regression [\citet{Efretal04}],
while it does not provide competitive performance relative to the
classical least-squares method with the DP, even when the true
number of change-points is known. To improve the performance, the
so-called LSTV* was further developed by incorporating a reduced
version of the DP. Roughly speaking, the LSTV plays essentially a
similar role in the LSTV* as our screening procedure in the NMCD. We
conduct comparisons between LSTV, LSTV* and NMCD under model (I)
only as the former two methods are not effective for scale changes
in model (II). The LSTV procedure is implemented until the
cardinality of the active set is exactly $K_n=11$, and both the NMCD
and LSTV* procedures are implemented by setting $\widebar{K}_{n}=30$ and
using the BIC to estimate the number of change-points.

%
\begin{table}
\tabcolsep=2pt
\caption{Comparison of NMCD, LSTV and LSTV* under model \textup{(I)}}\label{TLSTV}
\fontsize{7pt}{10pt}\selectfont
\begin{tabular*}{\tablewidth}{@{\extracolsep{\fill}}@{}ld{1.2}cccccccc@{}}
\hline
& & \multicolumn{2}{c}{\textbf{LSTV}} &\multicolumn{3}{c}{\textbf{LSTV*}} & \multicolumn{3}{c@{}}{\textbf{NMCD}}\\[-6pt]
& & \multicolumn{2}{c}{\hrulefill} &\multicolumn{3}{c}{\hrulefill} & \multicolumn{3}{c@{}}{\hrulefill}
\\
$\bolds{n}$& \multicolumn{1}{c}{$\bolds{\sigma}$} &
$\bolds{\xi(\widehat{\mathcal{G}}_n\|\mathcal{C}_t)}$&
$\bolds{\xi(\mathcal{C}_t\|\widehat{\mathcal{G}}_n)}$&
$\bolds{\xi(\widehat{\mathcal{G}}_n\|\mathcal{C}_t)}$&
$\bolds{\xi(\mathcal{C}_t\|\widehat{\mathcal{G}}_n)}$&
$\bolds{|\widehat{K}_n\,{-}\,K_n|}$&
$\bolds{\xi(\widehat{\mathcal{G}}_n\|\mathcal{C}_t)}$&
$\bolds{\xi(\mathcal{C}_t\|\widehat{\mathcal{G}}_n)}$&
$\bolds{|\widehat{K}_n\,{-}\,K_n|}$\\
\hline
\phantom{0}500 &  0.1  & 20.2 & 31.2 &  1.14 & 1.88 & 0.18 & 0.00 & 0.00 & 0.00\\
&  0.25  & 23.4 & 29.4 &  1.08 & 2.05 & 0.18 & 0.07 & 0.07 & 0.00 \\
&  0.5  & 26.1 & 27.0 &  2.10 & 3.14 & 0.17 & 1.39 & 1.30 & 0.03 \\
1000 &  0.1  & 43.1 & 60.2  & 2.82 & 2.21 & 0.15 & 0.00 & 0.00 & 0.00\\
&  0.25  & 46.2 & 59.4 &  3.23 & 2.24 & 0.16 & 0.04 & 0.04 & 0.00 \\
&  0.5  & 48.4 & 51.0 &  4.45 & 2.49 & 0.17 & 1.20 & 1.20 & 0.01
\\[3pt]
\multicolumn{4}{@{}l}{Computation time per run} &
\multicolumn{3}{c}{$n=500\dvtx 0.102$} &
\multicolumn{3}{c}{$n=500\dvtx 0.054$}
\\
\multicolumn{4}{@{}l}{(in seconds)} & \multicolumn{3}{c}{$n=1000\dvtx 0.776$} & \multicolumn{3}{c@{}}{$n=1000\dvtx 0.24$}\\
\hline
\end{tabular*}
\normalsize
\end{table}

The\vspace*{2pt} results in Table~\ref{TLSTV} show that the proposed NMCD and
LSTV* substantially outperform LSTV in terms of both
$\xi(\widehat{\mathcal{G}}_n\|\mathcal{C}_t)$ and
$\xi(\mathcal{C}_t\|\widehat{\mathcal{G}}_n)$. Moreover, the NMCD performs
uniformly better than LSTV*, which may be partly explained by the
fact that the induced shrinkage of LASSO often results in
significant bias toward zero for large regression coefficients
[\citet{FanLi01}]. Consequently, the LSTV also suffers from such bias,
which in turn may lead to unsatisfactory estimation of the locations
$\tau_k$'s. In Table~\ref{TLSTV}, we also report the average
computation time of the NMCD and LSTV* methods
using an Intel Core 2.2~MHz CPU. For a large sample size, NMCD is
much faster.

\subsection{Performance of NMCD with a diverging number of
change-points}\label{sec4.6}
To examine the setting that the number of
change-points increases with the sample size,
we choose seven increasing sample sizes,
$n={}$1000, 1500, 2000, 3000, 5000, 7500 and 10,000,
under models (I) and (II), respectively.
The number of change-points in model (I) is chosen as
$K_n=\lceil0.4n^{1/2}\rceil$, corresponding to the
values of 13, 16, 18, 22, 29, 35 and 40.
In each replication, we randomly generate the jump sizes $h_j$ as
follows: $h_{2k-1}=-1.5+\nu_{2k-1}$ and $h_{2k}=1.5+\nu_{2k}$,
$k=1,\ldots,\lceil K_n/2\rceil$, where $\nu_{j}\sim N(0,0.2^2)$. In
model (II), we take $K_n=\lceil0.2n^{1/2}\rceil$, and we only
consider the
scale changes (i.e., $h_j=0$ for all $j$) and the inflation (deflation)
sizes $v_j$ are chosen as: $v_{2k-1}=1/(5+\nu_{2k-1})$ and
$v_{2k}=5+\nu_{2k}$, $k=1,\ldots,\lceil K_n/2\rceil$,\vspace*{1pt} where
$\nu_{j}\sim N(0,0.2^2)$. We take the error distributions
to be $t_{(3)}$ and $\chi^2_{(1)}$ in models (I)~and~(II),
respectively. We fix $\sigma=0.5$, and generate $\{\tau_j/n\}_{j=1}^{K_n}$
from $U(0,1)$. All the tuning parameters are
the same as those in Section~\ref{sec4.3}.

%
\begin{figure}

\includegraphics{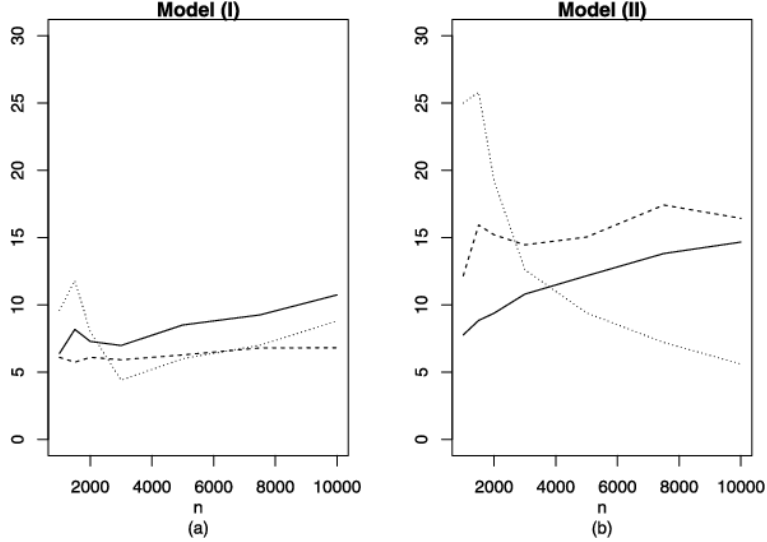}

\caption{The performance of NMCD under models
\textup{(I)} and \textup{(II)} when the number of change-points increases with the
sample size: the solid, dashed and dotted lines represent
$\xi(\widehat{\mathcal{G}}_n\|\mathcal{C}_t)$,
$\xi(\mathcal{C}_t\|\widehat{\mathcal{G}}_n)$, and $100|\widehat{K}_n-K_n|$
versus the sample size, respectively.}\label{fig-diverg}
\end{figure}

Figure~\ref{fig-diverg} depicts
the curves of $\xi(\widehat{\mathcal{G}}_n\|\mathcal{C}_t)$,
$\xi(\mathcal{C}_t\|\widehat{\mathcal{G}}_n)$, and $100|\widehat{K}_n-K_n|$
versus the sample size, respectively. For both models, all the
distance values are reasonably small and the three curves are
generally stable. This demonstrates that the NMCD is able to deliver
satisfactory detection performance with a diverging number of
change-points. From all these numerical studies, we conclude that
the proposed NMCD is a viable alternative approach to the MCP if
we take into account its efficiency, computational speed, and
robustness to error distributions and change patterns.

\section{Example}\label{sec5}

\begin{figure}

\includegraphics{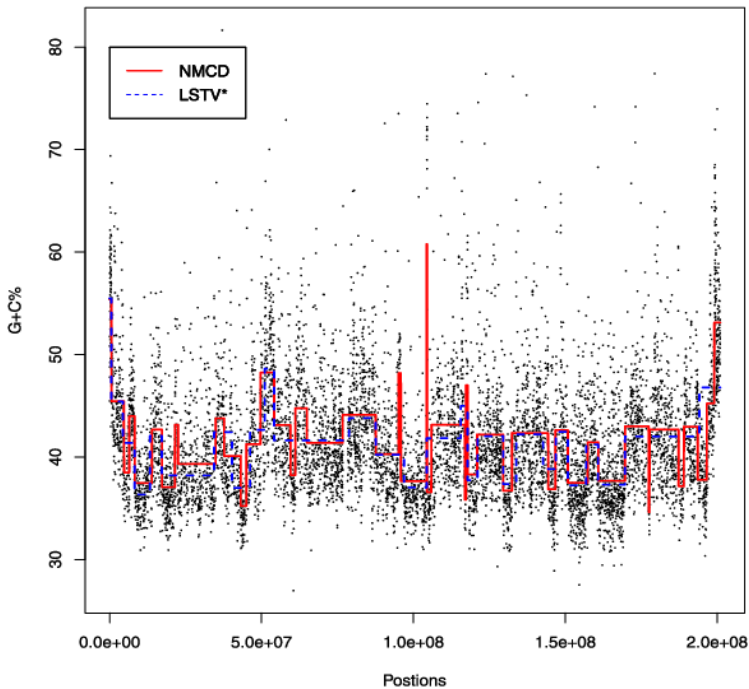}

\caption{Illustration of a chromosome
sequence with long homogeneous genome regions characterized by the
mean G${}+{}$C contents, together with the estimated changepoints using
the proposed NMCD and LSTV*, respectively. The red and blue solid
lines represent the sample means in each segmentation estimated by
NMCD and LSTV*, respectively.}\label{Isochore}
\end{figure}

For illustration, we apply the proposed NMCD procedure to identify
changes in the isochore structure, which refers to the proportion of
the G${}+{}$C composition in the large-scale DNA bases rather than A or T
[\citet{Olietal04}; \citet{FeaVas09}]. Such genetic
information is important to understand the evolution of base
composition, mutation and recombination rates. Figure~\ref{Isochore}
shows the G${}+{}$C content in percentage of a chromosome sequence with
long homogeneous genome regions characterized by well-defined mean
G${}+{}$C contents.

%

%
\begin{figure}

\includegraphics{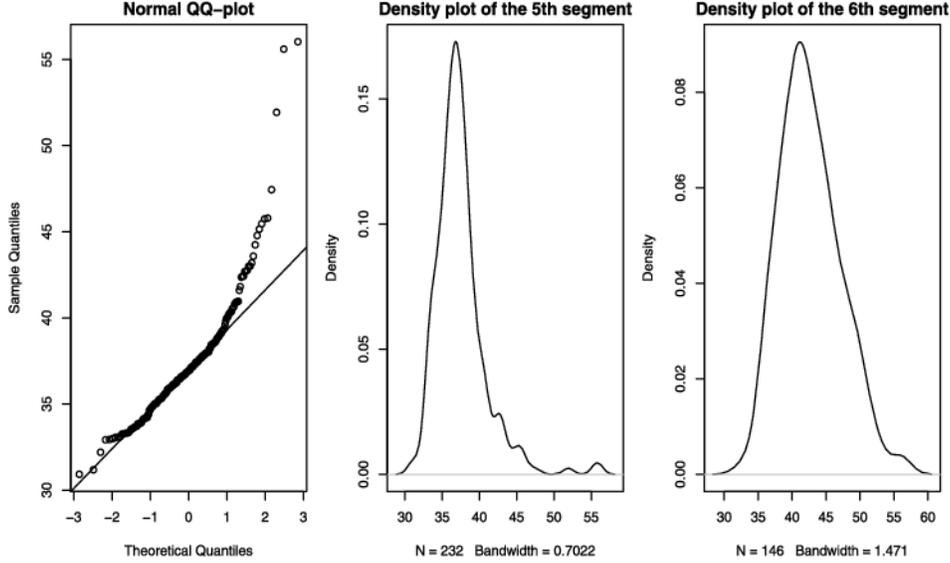}

\caption{The first plot: normal QQ-plot of
the 5th segment by using the NMCD; the second plot: density
estimation of the 5th segment; the third plot: density estimation of
the 6th segment.}\label{Figegana}
\end{figure}

As the data sequence appears to be complicated without any obvious
pattern and the sample size is large with $n=8811$, identification
of multiple change-points is very challenging. The data appear to
contain quite a few outlying observations, and thus we expect
that our nonparametric scheme would produce more robust detection results.

We take\vspace*{1pt} the upper bound for the number of change-points as
$\widebar{K}_{n}=100$, and set $n_{I}=\lceil(\log n)^{3/2}/2\rceil
=14$ and
$\zeta_n=(\log n)^{2}/2 \approx41$. After the initial screening
procedure, 305 candidate points remain, which dramatically reduces
the dimensionality of change-point detection. The BIC selection
criterion further leads to the estimated number of change-points
$\widehat{K}_n=43$. The entire procedure is completed in~54~seconds
using an Intel Core 2.2~MHz CPU. It can be seen from Figure~\ref
{Isochore} that the change-point estimates are generally
reasonable based on the proposed NMCD procedure. It can detect some
local and sharp features as well as those long unchanged data
segments. For comparison, we also apply the LSTV*
to the same dataset, and exhibit the result in Figure~\ref{Isochore}.
The estimated number of change-points using LSTV* is
$\widehat{K}_n=26$. We can see that both methods perform well, and the
line segments of
the two methods are largely overlapping, except that the NMCD tends
to detect relatively more picks or sharp changes.
Some large changes could be overlooked by LSTV* due to the LASSO-type
bias for large coefficients. This also explains
that the number of change-points identified by the LSTV* is smaller
than that of the proposed NMCD.

We performed the Shapiro--Wilk goodness-of-fit tests for normality
on the 44 segments identified by NMCD and found that 34 tests are
significant under the 0.01 nominal level. As an example, Figure~\ref
{Figegana} shows the normal QQ-plot of the fifth segment, from
which we can conclude that its distribution is far from normal.
Furthermore, the density estimation of two consecutive segments (the
5th and 6th) shown in Figure~\ref{Figegana} indicates that the two
distributions differ not only in the location but also in the scale
and shape. In light of these characteristics, our NMCD procedure is
more desirable than those parametric methods which need to specify
the mean or scale changes in advance.

\section{Concluding remarks}\label{sec6}

In the MCP, we have proposed a
nonparametric likelihood-based method for detection of multiple change-points.
The consistency of the proposed NMCD procedure is established under
mild conditions. The true number of change-points is assumed to be
unknown, and the BIC is used to choose the number of change-points.
To facilitate the implementation of NMCD, we suggest a DP
algorithm in conjunction with a screening procedure, which has been
shown to work well, particularly in large datasets. The computational
scheme is fast and competitive
with existing methods and, furthermore,
numerical comparisons show that
NMCD is able to strike a better balance for over- and
under-segmentation errors with nonnormal data and even has
comparable performance with the parametric model under the correctly
specified distributional assumption.

The proposed method is based on the assumption that there exists at
least one change point. In practical applications, we need to use
some tests within the nonparametric context to verify this
assumption. The tests proposed by \citet{EinMcK03} and \citet
{Zouetal07} are suited for this purpose. Our proposed NMCD is an
omnibus method,
and thus cannot diagnose whether a change occurs in the location,
scale, or shape.
To further determine which parameter changes, additional
nonparametric tests need to be used as an auxiliary tool. Moreover,
research is warranted to extend our method to other settings, such
as the autocorrelated observations, multivariate cases [\citet{MatJam}], and multiple structural changes in linear models [\citet{BaiPer98}].

\begin{appendix}
\section*{Appendix}\label{app}
First of all, we present a lemma in Wellner (\citeyear{We78}).
Let $G_n(u)$ denote the empirical C.D.F. of a random sample of $n$
uniform random variables on $(0,1)$, and define
$\|G_n(u)/u\|_s^t\equiv\sup_{s\leq u\leq t}(G_n(u)/u)$ and
\mbox{$G_n^{-1}(u)=\inf\{s\dvtx G_n(s)\geq u\}$}.

%
\begin{lemma}\label{wl}
For all $\lambda\geq0$ and $0\leq a\leq1$,
\begin{longlist}[(iii)]
\item[(i)] $\Pr(\|G_n(u)/u\|_a^1\geq\lambda)\leq
\exp\{-nah(\lambda)\}$,\vspace*{2pt}
\item[(ii)] $\Pr(\|u/G_n(u)\|_a^1\geq\lambda)\leq
\exp\{-nah(1/\lambda)\}$,\vspace*{2pt}
\item[(iii)] $\Pr(\|u/G_n^{-1}(u)\|_a^1\geq\lambda)\leq
\exp\{-naf(1/\lambda)\}$,\vspace*{2pt}
\item[(iv)] $\Pr(\|G_n^{-1}(u)/u\|_a^1\geq\lambda)\leq
\exp\{-naf(\lambda)\}$,\vspace*{2pt}
\item[(v)] $\Pr(\llvert\|G_n(u)/u-1\|_a^1\rrvert\geq\lambda
)\leq
2\exp(-nah(1+\lambda))$,
\end{longlist}
where $h(x)=x(\log x-1)+1$ and $f(x)=x+\log(1/x)-1$.
\end{lemma}

Before proceeding further, we state a key lemma, which allows us to
control the supremum of the likelihood function.

%
\begin{lemma}\label{klem4}
Suppose that assumptions \textup{(A1)--(A2)} hold and $K_n(\log\delta
_n)/\break \delta_n\rightarrow0$. Let
$w_n\equiv C_{\epsilon}K_n(\log K_n)^2(\log(\delta_nK_n))^2$, then
\[
\lim_{n\rightarrow\infty}K_n\Pr\Bigl\{\sup_{\tau_{m-1}\leq
k<l<\tau_{m-1}+\delta_n}
\xi_m(k,l)\geq w_n \Bigr\}<\epsilon,
\]
where
\begin{eqnarray*}
\xi_m(k,l)&=&n_{kl}\int_{X_{(1)}}^{X_{(n)}}
\biggl\{{\widehat{F}}_{k}^l(u)\log\biggl(
\frac{{\widehat
{F}}_{k}^l(u)}{F_m(u)} \biggr)
\\
&&\hspace*{45pt}{}+\bigl(1-{\widehat{F}}_{k}^l(u)\bigr) \log\biggl(
\frac{1-{\widehat{F}}_{k}^l(u)}{1-F_m(u)} \biggr) \biggr\} \frac
{d{\widehat{F}}
_n(u)}{{\widehat{F}}_n(u)(1-{\widehat{F}}_n(u))},
\end{eqnarray*}
$n_{kl}=l-k$ and $C_{\epsilon}$ is given in the proof.
\end{lemma}

\begin{pf}
Without loss of generality, suppose that $F_m$ is uniform on $[0,1]$
and $0<X_1<\cdots<X_n<1$. Then we have
%
%
\begin{equation}
\label{zct} \xi_m(k,l)=n_{kl}\int_{X_{(1)}}^{X_{(n)}}H
\bigl({\widehat{F}}_k^l(u),u\bigr)\bigl\{ {\widehat{F}}
_n(u) \bigl(1-{\widehat{F}}_n(u)\bigr)\bigr
\}^{-1}\,d{\widehat{F}}_n(u),
\end{equation}
where
\[
H(x,y)= x\log\biggl(\frac{x}{y} \biggr)+(1-x)\log\biggl(
\frac
{1-x}{1-y} \biggr).
\]
By setting $a_n=3h^{-1}(1+\alpha)\delta_n^{-1}\log(\delta
_nK_n)\equiv
D_{\alpha}\delta_n^{-1}\log(\delta_nK_n)$, $0<\alpha<1/2$, and
noting that
$h(1+\alpha)>0$, we write
\begin{eqnarray*}
&& \xi_m(k,l)
\\
&&\quad = n_{kl} \biggl(\int_{X_{(1)}}^{a_n}+
\int^{1-a_n}_{a_n}+\int_{1-a_n}^{X_{(n)}}
\biggr) H\bigl({\widehat{F}}_k^l(u),u\bigr)\bigl\{{
\widehat{F}}_n(u) \bigl(1-{\widehat{F}}_n(u)\bigr)\bigr\}
^{-1}\,d{\widehat{F}}_n(u)
\\
&&\quad \equiv \Delta_1+\Delta_2+\Delta_3.
\end{eqnarray*}

First, we provide an upper bound for $K_n\Pr(\sup_{k,l}\Delta_1\geq
w_n/3)$, where $\Delta_1\equiv\Delta_{11}+\Delta_{12}$ with
\begin{eqnarray*}
\Delta_{11}&=&n_{kl}\int_{X_{(1)}}^{a_n}
\frac{{\widehat
{F}}_{k}^l(u)}{u}\log\biggl(\frac{{\widehat{F}}_{k}^l(u)}{u} \biggr
)\frac{u}{{\widehat
{F}}_n(u)(1-{\widehat{F}}
_n(u))}\,d{
\widehat{F}}_n(u),
\\
\Delta_{12}&=&n_{kl}\int_{X_{(1)}}^{a_n}
\frac{1-{\widehat{F}}
_{k}^l(u)}{1-u}\log\biggl(\frac{1-{\widehat{F}}_{k}^l(u)}{1-u} \biggr
)\frac
{1-u}{{\widehat{F}}_n(u)(1-{\widehat{F}}_n(u))}\,d{
\widehat{F}}_n(u).
\end{eqnarray*}\eject\noindent
To show this, we choose $\lambda_{\epsilon1}$ such that as
$n\rightarrow\infty$,
\begin{eqnarray*}
&& K_n\Pr\bigl(\bigl\|u/{\widehat{F}}_n(u)
\bigr\|_{X_{(1)}}^1>\log K_n\lambda_{\epsilon1} \bigr)
\\
&&\qquad \leq K_n\Pr\biggl(\biggl\llVert\frac{nu}{(\tau_m-\tau_{m-1}){\widehat
{F}}_{\tau
_{m-1}}^{\tau_{m}}(u)}\biggr\rrVert
_{X_{(1)}}^1>\lambda_{\epsilon1}\log K_n \biggr)
\\
&&\qquad \leq n^{-1}(\tau_m-\tau_{m-1})K_n^2e
\lambda_{\epsilon1}\exp\bigl\{ -n^{-1}(\tau_m-
\tau_{m-1})\lambda_{\epsilon1} \log K_n\bigr\}<
\epsilon/12,
\end{eqnarray*}
based on assumption (A2) and the fact that
\[
\Pr\bigl(\bigl\|u/G_n(u)\bigr\|_{X_{(1)}}^1>\lambda\bigr)\leq
\Pr\bigl(\bigl\|G_n^{-1}(u)/u\bigr\| _{1/n}^1
\geq\lambda\bigr)\leq e\lambda\exp\{-\lambda\}
\]
by using Lemma \ref{wl}(iv). Similarly,
\[
K_n\Pr\bigl\{\bigl\|{\widehat{F}}_n^{-1}(u)/u
\bigr\|_{1/n}^1>\lambda_{\epsilon
1}\log K_n \bigr
\}<\epsilon/12.
\]

Also, we consider the event
$A_m\equiv\bigcup_{k,l} \{\|{\widehat{F}}_k^l(u)/u\|_0^1>
\lambda_{\epsilon2}K_n{\delta_n}/{n_{kl}} \}$, and thus
\begin{eqnarray*}
K_n\Pr(A_m)&=&K_n\Pr\biggl(\bigcup
_{k,l}\frac{n_{kl}}{\delta_n}\bigl\| {\widehat{F}} _k^l(u)/u
\bigr\|_0^1>\lambda_{\epsilon2}K_n \biggr)
\\
&\leq& K_n\Pr\biggl(\bigcup_{k,l}\bigl\|{
\widehat{F}}_{\tau_{m-1}}^{\tau
_{m-1}+\delta_n}(u)/u\bigr\|_0^1>
\lambda_{\epsilon2}K_n \biggr)
\\
&=&K_n\Pr\bigl(\bigl\|{\widehat{F}}_{\tau_{m-1}}^{\tau_{m-1}+\delta
_n}(u)/u\bigr\|
_0^1>\lambda_{\epsilon2}K_n \bigr) \leq
e\lambda_{\epsilon2}^{-1}<\epsilon/12
\end{eqnarray*}
by choosing a proper $\lambda_{\epsilon2}$. In parallel,
let $B_m\equiv\bigcup_{k,l} \{\|(1-{\widehat{F}}_k^l(u))/(1-u)\|_0^1>
\lambda_{\epsilon2}{K_n\delta_n}/{n_{kl}} \}$, and we have
\begin{eqnarray*}
&&K_n\Pr\bigl(\bigl\|(1-u)/\bigl(1-{\widehat{F}}_n(u)\bigr)
\bigr\|_0^1>\lambda_{\epsilon
1}\log K_n
\bigr)<\epsilon/12
\end{eqnarray*}
and $K_n\Pr(B_m)<e\lambda_{\epsilon2}^{-1}<\epsilon/12$.

For the interaction of the events $\widebar{A}_m$,
$\|u/{\widehat{F}}_n(u)\|_0^1\leq\lambda_{\epsilon1}\log K_n$, and
\[
\bigl\|{\widehat{F}}_n^{-1}(u)/u\bigr\|_{1/n}^1\leq\lambda_{\epsilon1}\log
K_n,
\]
we have
\begin{eqnarray*}
\Delta_{11}&=&n_{kl}\int_{X_{(1)}}^{a_n}
\frac{{\widehat
{F}}_{k}^l(u)}{u}\log\biggl(\frac{{\widehat{F}}_{k}^l(u)}{u} \biggr
)\frac{u}{{\widehat
{F}}_n(u)}
\frac
{1}{(1-{\widehat{F}}_n(u))}\,d{\widehat{F}}_n(u)
\\
&\leq& -n_{kl}\frac{K_n\delta_n}{n_{kl}}\lambda_{\epsilon2}\log\biggl(
\frac{K_n\delta_n}{n_{kl}}\lambda_{\epsilon2} \biggr) \lambda_{\epsilon
1}\log
K_n\log\bigl(1-{\widehat{F}}_n^{-1}(a_n)
\bigr)
\\
&\leq&-{K_n\delta_n}\lambda_{\epsilon2}\log
({K_n\delta_n}\lambda_{\epsilon2} )
\lambda_{\epsilon1}\log K_n \log(1-\lambda_{\epsilon1}a_n
\log K_n)
\\
&\leq& K_n(\log K_n)^2\delta_na_n
\lambda_{\epsilon2}\lambda_{\epsilon_1}^2\log(
\delta_nK_n) \bigl(1+o(1)\bigr)
\end{eqnarray*}
as $n\rightarrow\infty$. Consequently, as $n\rightarrow\infty$,
\begin{eqnarray*}
&& K_n\Pr\Bigl(\sup_{k,l}\Delta_{11}
\geq w_n/6 \Bigr)
\\
&&\qquad \leq K_n\Pr(A_m)+K_n\Pr\bigl(\bigl\|{
\widehat{F}}_n^{-1}(u)/u\bigr\| _{1/n}^1>
\lambda_{\epsilon1}\log K_n\bigr)
\\
&&\quad\qquad{} +K_n\Pr\bigl(\bigl\|u/{\widehat{F}}_n(u)
\bigr\|_{X_{(1)}}^1>\lambda_{\epsilon
1}\log K_n\bigr)
+\delta_n^2K_n\Pr(\Delta_{11}\geq
w_n/6)
\\
&&\qquad \leq \frac{1}{4}\epsilon+\delta_n^2K_n
\Pr\bigl\{\lambda_{\epsilon2}\lambda_{\epsilon_1}^2(\log
\delta_nK_n)^2K_n(\log
K_n)^2\bigl(1+o(1)\bigr)\geq w_n/6 \bigr\}
= \frac{1}4\epsilon,
\end{eqnarray*}
where the probability $\Pr\{(\log
\delta_nK_n)^2K_n(\log K_n)^2\lambda_{\epsilon2}\lambda_{\epsilon
1}^2(1+o(1))\geq w_n \}$
would be zero when $n$ is
sufficiently large, as long as
$C_{\epsilon}>6D_{\alpha}\lambda_{\epsilon2}\lambda_{\epsilon_1}^2$.\vspace*{1pt}

Similarly, we can show that $K_n\Pr(\sup_{k,l}\Delta_{12}\geq
w_n/6)\leq
\epsilon/4$ as $n\rightarrow\infty$. Thus,
\[
\Pr\Bigl(\sup_{k,l}\Delta_{1}\geq
w_n/3 \Bigr)\le\Pr\Bigl(\sup_{k,l}
\Delta_{11}\geq w_n/6 \Bigr) +\Pr\Bigl(\sup
_{k,l}\Delta_{12}\geq w_n/6 \Bigr)<
\epsilon/2.
\]
By symmetry, we immediately have
\[
K_n\Pr\Bigl(\sup_{k,l}\Delta_{3}\geq
w_n/3 \Bigr)\leq\frac{1}2\epsilon\qquad\mbox{as } n\rightarrow
\infty.
\]
Thus, it remains to give a bound of $K_n\Pr(\sup_{k,l}\Delta_2\geq
w_n/3)$.
Following similar argument in the proof of Theorem~3.1 of \citet
{JagWel07}, we can express $H({\widehat{F}}_{kl}(u),u)$ as
\[
H\bigl({\widehat{F}}_n(u),u\bigr)=\frac{1}2
\frac{({\widehat
{F}}_{k}^l(u)-u)^2}{{\widehat{F}}_{kl}^*(u)(1-{\widehat{F}}
_{kl}^*(u))}
\]
for $0<u<1$ where $|{\widehat{F}}_{kl}^*(u)-u|\leq|{\widehat
{F}}_{k}^l(u)-u|$. Then we
rewrite $\Delta_2$ as
\begin{eqnarray*}
\Delta_2&=&\frac{1}2\int_{a_n}^{1-a_n}
\frac{n_{kl}({\widehat
{F}}_{k}^l(u)-u)^2}{u(1-u)} \frac{u(1-u)}{{\widehat
{F}}_{kl}^*(u)(1-{\widehat
{F}}_{kl}^*(u))}\frac{d{\widehat{F}}_n(u)}{{\widehat{F}}
_n(u)(1-{\widehat{F}}_n(u))}
\\
&\leq&\frac{1}2\biggl\llVert\frac{n_{kl}({\widehat
{F}}_{k}^l(u)-u)^2}{u(1-u)}\biggr\rrVert
_{a_n}^{1-a_n} \biggl\llVert\frac{u}{{\widehat{F}}_{kl}^*(u)}\biggr
\rrVert
_{a_n}^{1-a_n}\biggl\llVert\frac
{1-u}{1-{\widehat{F}}_{kl}^*(u)}\biggr\rrVert
_{a_n}^{1-a_n}
\\
&&{}\times  \int_{a_n}^{1-a_n}
\frac{d{\widehat{F}}_n(u)}{{\widehat
{F}}_n(u)(1-{\widehat{F}}_n(u))}.
\end{eqnarray*}
Consider the event
$C_m\equiv\bigcup_{k,l} \{ |\|{\widehat{F}}_k^l(u)/u-1\|
_{a_n}^{1-a_n} |>\alpha\}$
for some $0<\alpha<1$ and, by applying Lemma \ref{wl}(v), we have
\begin{eqnarray*}
K_n\Pr(C_m)&\leq&\delta_n^2K_n
\Pr\bigl( \bigl|\bigl\|{\widehat{F}}_{k}^l(u)/u-1\bigr\|
_{a_{n}}^{1-a_n} \bigr|>\alpha\bigr)
\\
&\leq&2\exp\bigl\{2\log(\delta_nK_n)-
\delta_na_nh(1+\alpha)\bigr\} \rightarrow0.
\end{eqnarray*}
On the event $\widebar{C}_m$ and
$|{\widehat{F}}_{kl}^*(u)/u-1|<|{\widehat{F}}_{kl}(u)/u-1|<\alpha$,
we have
\[
\biggl\llVert\frac{u}{{\widehat{F}}_{kl}^*(u)}\biggr\rrVert_{a_n}^{1-a_n}<
\frac
{1}{1-\alpha}.
\]
Symmetrically, we also have
\[
\biggl\llVert\frac{1-u}{1-{\widehat{F}}_{kl}^*(u)}\biggr\rrVert_{a_n}^{1-a_n}<
\frac
{1}{1-\alpha}
\]
on the event $\widebar{D}_m$, where
$D_m\equiv\bigcup_{k,l} \{ |\|(1-{\widehat
{F}}_k^l(u))/(1-u)-1\|
_{a_n}^{1-a_n} |>\alpha\}$
occurs with the probability tending to zero. On the other hand, by
using Lemma \ref{wl}(v) again, it is easy to see that, for sufficiently large $n$,
\[
\int_{a_n}^{1-a_n}\bigl\{{\widehat{F}}_n(u)
\bigl(1-{\widehat{F}}_n(u)\bigr)\bigr\} ^{-1}\,d{
\widehat{F}}_n(u)\leq-2\log a_n+C_\alpha\le2\log(
\delta_nK_n) \bigl(1+o_p(1)\bigr),
\]
where the constant $C_\alpha$ depends on $\alpha$.

Now, we consider the term
$\|{n_{kl}({\widehat{F}}_{k}^l(u)-u)^2}/\{u(1-u)\}\|_{a_n}^{1-a_n}$,
and let
$\varrho_n=(w_n/\log(\delta_nK_n))^{1/2}$. By taking $q(t)=\sqrt
{t(1-t)}$ in
Inequality 11.2.1 of \citeauthor{ShoWel86} [(\citeyear{ShoWel86}), page~446],
\begin{eqnarray*}
\Pr\biggl(\biggl\llVert\frac{n_{kl}({\widehat{F}}_{k}^l(u)-u)^{\pm
}}{\sqrt{u(1-u)}}\biggr\rrVert_{a_n}^{1/2}
\geq\varrho_n \biggr) &\leq&6\int_{a_n}^{1/2}
\frac{1}t\exp\biggl\{-\frac{1}8\gamma^{\pm}\varrho
_n^2(1-t) \biggr\}\,dt
\\
&\leq&6\exp\biggl\{-\frac{1}{16}\gamma^{\pm}\varrho_n^2
\biggr\} \log\delta_n\bigl(1+o(1)\bigr),
\end{eqnarray*}
where $\gamma^{-}=1$, $\gamma^{+}=\psi(\varrho_n/\sqrt{\delta
_na_n})$, and
$\psi(x)=2h(1+x)/x^2$. By using the fact that $\psi(x)\sim2(\log
x)/x$ as $x\rightarrow\infty$ [Proposition 11.1.1 in \citet{ShoWel86}],
$\gamma^{+}\sim\log(C_{\epsilon}K_n(\log
K_n)^2)/(C_{\epsilon}^{1/2}K_n^{1/2}\log K_n)$ for
sufficiently large $C_{\epsilon}$. Consequently, we have
\begin{eqnarray*}
&& K_n\Pr \biggl(\sup_{k,l}\biggl\llVert
\frac{n_{kl}({\widehat{F}}
_{k}^l(u)-u)^2}{u(1-u)}\biggr\rrVert_{a_n}^{1/2}\geq
\frac{(1-\alpha)^2w_n}{3\log(\delta_nK_n)} \biggr)
\\
&&\qquad \leq K_n\delta_n^2\Pr\biggl(\biggl\llVert
\frac{n_{kl}({\widehat
{F}}_{k}^l(u)-u)^{\pm
}}{\sqrt{u(1-u)}}\biggr\rrVert_{a_n}^{1/2}\geq
\varrho_n \biggr)
\\
&&\qquad \leq 12\exp\biggl(2\log(\delta_nK_n)-\frac{1}{16}
\gamma^{+}\varrho_n^2 \biggr)\log
\delta_n\bigl(1+o(1)\bigr)
\\
&&\qquad \rightarrow0\qquad\mbox{as } \delta_n\rightarrow\infty
\end{eqnarray*}
as long as $C_{\epsilon}$ is sufficiently large. By symmetry, we can
also show that
\begin{eqnarray*}
&& K_n\Pr \biggl(\sup_{k,l}\biggl\llVert
\frac{n_{kl}({\widehat{F}}
_{k}^l(u)-u)^2}{u(1-u)}\biggr\rrVert_{1/2}^{1-a_n}\geq
\frac{(1-\alpha)^2w_n}{3\log(\delta_nK_n)} \biggr)\rightarrow0.
\end{eqnarray*}
Finally, we obtain as $n\rightarrow\infty$,
\begin{eqnarray*}
&& K_n\Pr\Bigl(\sup_{k,l}\Delta_2\ge
w_n/3 \Bigr)
\\
&&\quad \leq  K_n \Pr(C_m)+K_n
\Pr(D_m)
\\
&&\qquad{}+K_n\Pr\biggl(\sup_{k,l}\biggl\llVert
\frac{n_{kl}({\widehat{F}}
_{k}^l(u)-u)^2}{u(1-u)}\biggr\rrVert_{a_n}^{1-a_n}\frac{1}{(1-\alpha
)^2}
\log(\delta_nK_n)\geq w_n/3 \biggr)
\rightarrow0,
\end{eqnarray*}
which completes the proof of this lemma.
\end{pf}

By Lemma \ref{klem4}, the next lemma follows immediately.
%
%
\begin{lemma}\label{klem}
Suppose that assumptions \textup{(A1)--(A2)} hold and $K_n(\log n)/n\rightarrow0$.
Then
\[
\lim_{n\rightarrow\infty}K_n\Pr\Bigl\{\sup_{\tau_{m-1}\leq
k<l<\tau_{m}}
\xi_m(k,l)\geq u_n \Bigr\}<\epsilon,
\]
where $u_{n}\equiv C_{\epsilon}K_n(\log K_n)^2(\log(nK_n))^2$ with
a sufficiently large $C_{\epsilon}$.
\end{lemma}

Let $\widetilde{O}_p(q_n; K_n)$ be a sequence of positive random
variables $Z_n$ if for any $\epsilon>0$,
\[
\lim_{n\to\infty} K_n\Pr(Z_n>C_{\epsilon}q_n)<
\epsilon,
\]
where $C_{\epsilon}$ is a constant depending only on $\epsilon$.

%
\begin{lemma}\label{yl2}
Suppose that assumptions \textup{(A1)--(A2)} hold. For any $L\geq1$ and
$\tau_s<\tau_1'<\cdots<\tau_L'<\tau_{s+1}$, as $n\rightarrow
\infty$,
\begin{eqnarray*}
0&\leq& R_n\bigl(\tau_s,\tau_1',
\ldots,\tau_L',\tau_{s+1}\bigr)-R_n(
\tau_s,\tau_{s+1})
\\
&=&\widetilde{O}_p\bigl(L^2K_n\bigl(
\log(K_nL)\bigr)^2\bigl(\log(nK_nL)
\bigr)^2; K_n\bigr).
\end{eqnarray*}
\end{lemma}

\begin{pf}
By noting that $H(x,y)$ is a convex
function, the left
inequality is obvious. Without loss of generality, we assume $L=1$,
and for $L>1$ the result follows by induction. By the fact that
$(\tau_1'-\tau_s){\widehat{F}}_{\tau_s}^{\tau_1'}(u)+(\tau
_{s+1}-\tau_1')
{\widehat{F}}_{\tau_1'}^{\tau_{s+1}}(u)=(\tau_{s+1}-\tau
_s){\widehat{F}}_{\tau
_s}^{\tau_{s+1}}(u)$,
\begin{eqnarray*}
R_n\bigl(\tau_s,\tau_1',
\tau_{s+1}\bigr)-R_n(\tau_s,
\tau_{s+1})&=&\xi_{s}\bigl(\tau_s,
\tau_1'\bigr)+\xi_{s}\bigl(
\tau_1',\tau_{s+1}\bigr)-\xi_{s}(
\tau_s,\tau_{s+1})
\\
&\le&\xi_{s}\bigl(\tau_s,\tau_1'
\bigr)+\xi_{s}\bigl(\tau_1',
\tau_{s+1}\bigr).
\end{eqnarray*}
Similarly, for any $L$, we have
\begin{eqnarray*}
&& R_n\bigl(\tau_s,\tau_1',
\ldots,\tau_L',\tau_{s+1}\bigr)-R_n(
\tau_s,\tau_{s+1})
\\
&&\qquad  \le\xi_{s}\bigl(
\tau_s,\tau_1'\bigr)+\xi_{s}
\bigl(\tau_1',\tau_2'\bigr)+
\cdots+\xi_s\bigl(\tau_L',
\tau_{s+1}\bigr).
\end{eqnarray*}
Thus, for any $\epsilon>0$,
\begin{eqnarray*}
\hspace*{-3pt}&& \lim_{n\to\infty} K_n \Pr\bigl\{R_n\bigl(
\tau_s,\tau_1',\ldots,\tau
_L',\tau_{s+1}\bigr)-R_n(
\tau_s,\tau_{s+1})
\\
\hspace*{-3pt}&&\hspace*{54pt}{} >C_{\epsilon}L^2K_n
\bigl(\log(K_nL)\bigr)^2\bigl(\log(nK_nL)
\bigr)^2 \bigr\}
\\
\hspace*{-3pt}&&\qquad \le \lim_{n\to\infty} K_n \Pr\bigl\{\xi_{s}
\bigl(\tau_s,\tau_1'\bigr)+\cdots+\xi_s\bigl(\tau_L',\tau_{s+1}
\bigr)
\\
\hspace*{-3pt}&&\hspace*{84pt}{}  >C_{\epsilon}L^2K_n\bigl(\log(K_nL)
\bigr)^2\bigl(\log(nK_nL)\bigr)^2 \bigr\}
\\
\hspace*{-3pt}&&\qquad \le L^{-1}\sum_{k=0}^{L}\lim
_{n\to\infty} K_nL \Pr\bigl\{\xi_s\bigl(
\tau_k'-\tau_{k+1}'
\bigr)
\\
\hspace*{-3pt}&&\hspace*{96pt}\quad\qquad{} >C_{\epsilon}L^2K_n\bigl(\log(K_nL)
\bigr)^2\bigl(\log(nK_nL)\bigr)^2 \bigr\}
\\
\hspace*{-3pt}&&\qquad < L^{-1}(L+1) \epsilon,
\end{eqnarray*}
where the last result follows immediately from
Lemma \ref{klem}.
\end{pf}

Next, we demonstrate that the global minimum of the BIC includes no
less than
$K_n$ change-point estimators asymptotically.

%
\begin{pro}\label{prol}
If assumptions \textup{(A1)--(A4)} hold, $\Pr\{\widehat{K}_n\geq K_n\}
\rightarrow1$.
\end{pro}

\begin{pf}
Define $\rho_n=\lambda_n/8$, and consider $0<L<K_n$. Let
\begin{eqnarray*}
&& B_r(L,\rho_n)
\\
&&\qquad =\bigl\{\bigl(\tau_1',
\ldots,\tau_{L}'\bigr)\dvtx 1<\tau_1'<
\cdots<\tau_{L}'\leq n \mbox{ and } \bigl|\tau_s'-
\tau_r\bigr|>\rho_n\mbox{ for } 1\leq s\leq{L}\bigr\},
\end{eqnarray*}
$r=1,\ldots,K_n$. For $L<K_n$, $(\hat{\tau}_1,\ldots,\hat{\tau}_{L})$
must belong to one $B_r(L,\rho_n)$.
For every $(\tau_1',\ldots,\tau_{L}')\in B_r(L,\rho_n)$, we have
%
%
\begin{eqnarray}\label{yao10}
&& R_n\bigl(\tau_1',\ldots,\tau_{L}'\bigr)
\nonumber\\[-8pt]\\[-8pt]
&&\qquad \leq R_n\bigl(
\tau_1',\ldots,\tau_{L}',
\tau_1,\ldots,\tau_{r-1},\tau_r-
\rho_n, \tau_r+\rho_n,\tau_{r+1},
\ldots,\tau_{K_n}\bigr)\nonumber
\end{eqnarray}
and the right-hand side of (\ref{yao10}) can be expressed as
$T_1+\cdots+T_{K_n+2}$, where $T_s$ ($s=1,\ldots,r-1,r+2,\ldots,K_n+1$)
is the sum of integrals involving the $X_i$'s ($\tau_{s-1}\leq
i<\tau_s$); $T_r$ is that involving the $X_i$'s ($\tau_{r-1}\leq
i<\tau_r-\rho_n$); $T_{r+1}$ is that involving the
$X_i$'s ($\tau_r+\rho_n \leq i<\tau_{r+1}$); $T_{K_n+2}$
is that involving the $X_i$'s ($\tau_r-\rho_n \leq
i<\tau_r+\rho_n$). For
$s=1,\ldots,r-1,r+2,\ldots,K_n+1$, by Lemma \ref{yl2}, we have
\begin{eqnarray*}
R_n(\tau_{s-1},\tau_s)\leq T_s &\leq& R_n(\tau_{s-1},\tau_s)+
\widetilde{O}_p\bigl(L^2K_n\bigl(
\log(K_nL)\bigr)^2\bigl(\log(nK_nL)
\bigr)^2\bigr)
\\
&=&R_n(\tau_{s-1},\tau_s)+
\widetilde{O}_p(b_n;K_n),
\end{eqnarray*}
where $b_n=K_n^3(\log K_n)^2(\log n)^2$. Similarly, we have
\begin{eqnarray*}
T_r&=&R_n(\tau_{r-1},\tau_r-
\rho_n)+\widetilde{O}_p(b_n;K_n),
\\
T_{r+1}&=&R_n(\tau_r+\rho_n,
\tau_{r+1})+\widetilde{O}_p(b_n;K_n)
\end{eqnarray*}
and in addition,
\begin{eqnarray*}
T_{K_n+2}&=&R_n(\tau_r-\rho_n,{
\tau_r+\rho_n })+R_n(\tau_r+
\rho_n,\tau_{r+1})
\\
&=&R_n(\tau_r-\rho_n,\tau_r)+R_n(
\tau_r,\tau_r+\rho_n)+\Delta,
\end{eqnarray*}
where $\Delta\equiv R_n(\tau_r-\rho_n,{\tau_r+\rho_n})
-R_n(\tau_r-\rho_n,\tau_r)-R_n(\tau_r,\tau_r+\rho_n )$. Note
that
\begin{eqnarray*}
\Delta&=&2\rho_n \int_{X_{(1)}}^{X_{(n)}}\bigl[{
\widehat{F}}_{\tau_r-\rho
_n}^{\tau
_r+\rho_n}(u)\log\bigl(F_{r,1/2}(u)
\bigr)
\\
&&\hspace*{47pt}{} +\bigl\{1-{\widehat{F}}_{\tau_r-\rho
_n}^{\tau
_r+\rho_n}(u)\bigr\}\log
\bigl(1-F_{r,1/2}(u)\bigr)\bigr]\,dw(u)
\\
&&{}-\rho_n\int_{X_{(1)}}^{X_{(n)}}\bigl[{
\widehat{F}}_{\tau_r-\rho
_n}^{\tau
_r}(u)\log\bigl(F_{r}(u)\bigr)
\\
&&\hspace*{54pt}{} +
\bigl\{1-{\widehat{F}}_{\tau_r-\rho_n}^{\tau
_r}(u)\bigr\}\log
\bigl(1-F_{r}(u)\bigr)\bigr]\,dw(u)
\\
&&{}-\rho_n\int_{X_{(1)}}^{X_{(n)}}\bigl[{
\widehat{F}}_{\tau_r}^{\tau
_r+\rho
_n}(u)\log\bigl(F_{{r+1}}(u)\bigr)
\\
&&\hspace*{53pt}{}+
\bigl\{1-{\widehat{F}}_{\tau_r}^{\tau_r+\rho
_n}(u)\bigr\}\log
\bigl(1-F_{{r+1}}(u)\bigr)\bigr]\,dw(u)+\widetilde{O}_p(b_n;K_n)
\\
&=&-\rho_n\int_{X_{(1)}}^{X_{(n)}}\biggl[{
\widehat{F}}_{\tau_r-\rho
_n}^{\tau
_r}(u)\log\biggl(\frac{F_{r}(u)}{F_{r,1/2}(u)}
\biggr)
\\
&&\hspace*{53pt}{}+\bigl\{1-{\widehat{F}}_{\tau
_r-\rho_n}^{\tau_r}(u)\bigr\}\log
\biggl( \frac{1-F_{r}(u)}{1-F_{r,1/2}(u)} \biggr)\biggr]\,dw(u)
\\
&&{}-\rho_n\int_{X_{(1)}}^{X_{(n)}}\biggl[{
\widehat{F}}_{\tau_r}^{\tau
_r+\rho
_n}(u)\log\biggl(\frac{F_{\tau_{r+1}}(u)}{F_{r,1/2}(u)}
\biggr)
\\
&&\hspace*{54pt}{}+\bigl\{ 1-{\widehat{F}}_{\tau_r}^{\tau_r+\rho_n}(u)\bigr\}\log
\biggl(\frac
{1-F_{{r+1}}(u)}{1-F_{r,1/2}(u)} \biggr)\biggr]\,dw(u)
\\
&&{}+\widetilde{O}_p(b_n;K_n)
\\
&\equiv&-\widetilde{\Delta}+\widetilde{O}_p(b_n;K_n).
\end{eqnarray*}
Let $\widetilde{\Delta}=\widetilde{\Delta}_{1}+\widetilde{\Delta
}_2$, and then
\begin{eqnarray*}
\widetilde{\Delta}_{1}&\geq&\rho_n\int
_{X_{(1)}}^{X_{(n)}}\biggl[{\widehat{F}}_{\tau
_r-\rho_n}^{\tau_r}(u)
\log\biggl(\frac
{F_{r}(u)}{F_{r,1/2}(u)} \biggr)
\\
&&\hspace*{42pt}{} +\bigl\{1-{\widehat{F}}_{\tau_r-\rho
_n}^{\tau
_r}(u)
\bigr\}\log\biggl(\frac{1-F_{r}(u)}{1-F_{r,1/2}(u)} \biggr)\biggr]\,dw(u)
\\
&=&\rho_n\int_{0}^{1}\biggl[{
\widehat{F}}_{\tau_r-\rho_n}^{\tau_r}(u)\log\biggl(\frac
{F_{r}(u)}{F_{r,1/2}(u)}
\biggr)
\\
&&\hspace*{30pt}{}+\bigl\{1-{\widehat{F}}_{\tau
_r-\rho_n}^{\tau
_r}(u)\bigr\}\log
\biggl(\frac{1-F_{r}(u)}{1-F_{r,1/2}(u)} \biggr)\biggr]\,dw(u)
\\
&\equiv&\widetilde{\Delta}_{1}'.
\end{eqnarray*}
By assumption (A3), we have
\begin{eqnarray*}
\widetilde{\Delta}_{1}'&=&\rho_n\int
_{0}^{1} \biggl[F_{r}(u)\log\biggl(
\frac{F_{{r}}(u)}{F_{r,1/2}(u)} \biggr) +\bigl\{1-F_{r}(u)\bigr\}\log
\biggl(
\frac{1-F_{{r}}(u)}{1-F_{r,1/2}(u)} \biggr) \biggr]
\\
&&\hspace*{25pt}{}\times\frac{1}{F(u)(1-F(u))}\,dF(u) \bigl(1+o(1)\bigr),\qquad\mbox{a.s.}
\end{eqnarray*}
Using the similar procedure, we can obtain the corresponding bound
for $\widetilde{\Delta}_2$. As a result, as $n\to\infty$,
\begin{eqnarray*}
\widetilde{\Delta}&\geq& \rho_n \biggl\{\int_{0}^{1}
\biggl[{F_{r}}(u)\log\biggl(\frac{F_{{r}}(u)}{F_{r,1/2}(u)} \biggr)
+\bigl
\{1-F_{r}(u)\bigr\}\log\biggl(\frac{1-F_{{r}}(u)}{1-F_{r,1/2}(u)}
\biggr) \biggr]
\\
&&\hspace*{30pt}{}\times\frac{1}{F(u)(1-F(u))}\,dF(u)
\\
&&\hspace*{15pt}{}+\int_{0}^{1} \biggl[F_{{r+1}}(u)\log
\biggl(\frac
{F_{{r+1}}(u)}{F_{r,1/2}(u)} \biggr)+ \bigl\{1-F_{{r+1}}(u)\bigr\}\log
\biggl(
\frac
{1-F_{{r+1}}(u)}{1-F_{r,1/2}(u)} \biggr) \biggr]
\\
&&\hspace*{209pt}{}\times\frac{1}{F(u)(1-F(u))}\,dF(u) \biggr\}
\\
&\equiv&\rho_n S(F_{{r}},F_{{r+1}}),
\end{eqnarray*}
in which the distance $S(F_{{r}},F_{{r+1}})$ is strictly larger than zero.

Therefore,
\begin{eqnarray*}
&& \max_{(\tau_1',\ldots,\tau_{L}')\in
B_r(L,\rho_n)}R_n\bigl(\tau_1',
\ldots,\tau_{L}'\bigr)
\\
&&\qquad \leq \max_{(\tau_1',\ldots,\tau_{L}')\in
B_r(L,\rho_n)}R_n\bigl(\tau_1',
\ldots,\tau_{L}',\tau_1,\ldots,\tau
_{r-1},\tau_r-\rho_n,\tau_r+ \rho_n,
\\
&&\hspace*{248pt} \tau_{r+1},\ldots,\tau_{K_n}\bigr)
\\
&&\qquad =\sum_{s\neq
r,r+1}^{K_n+1}R_n(
\tau_{s-1},\tau_s)+R_n(\tau_{r-1},
\tau_r-\rho_n)+R_n(\tau_r-
\rho_n,\tau_r)
\\
&&\quad\qquad{}+R_n(\tau_r,
\tau_r+\rho_n)+R_n(\tau_r+\rho_n,\tau_{r+1})+
\Delta+\widetilde{O}_p(b_n;K_n)
\\
&&\qquad \leq R_n(\tau_1,\ldots,\tau_{K_n})-
\rho_n S(F_{{r}},F_{{r+1}})+\widetilde{O}_p(b_n;K_n).
\end{eqnarray*}
Let $\mathrm{BIC}_{*}=-R_n(\tau_1,\ldots,\tau_{K_n})+K_n\zeta_n$,
and for
$L<K_n$, with probability tending to 1, we have
\[
\mathrm{BIC}_L-\mathrm{BIC}_{*}\geq\rho_n
S(F_{{r}},F_{{r+1}})-\widetilde{O}_p(b_n;K_n)-(K_n-L)
\zeta_n
\]
as $n\rightarrow\infty$. For any $\epsilon>0$, we have, as $n\to
\infty$,
\begin{eqnarray*}
\Pr(\widehat{K}_n<K_n)&=&\Pr\Biggl(\bigcup
_{L=1}^{K_n-1}(\mathrm{BIC}_L<
\mathrm{BIC}_{*}) \Biggr)\le\sum_{L=1}^{K_n-1}
\Pr(\mathrm{BIC}_L<\mathrm{BIC}_{*} )
\\
&\le&\sum_{L=1}^{K_n-1}\Pr\bigl(
\widetilde{O}_p(b_n; K_n)>\rho
_nS(F_r,F_{r+1})-(K_n-L)
\zeta_n \bigr)
\\
&\le&K_n\Pr\bigl(\widetilde{O}_p(b_n;K_n)>b_n
\bigr)<\epsilon.
\end{eqnarray*}
This completes the proof of this proposition.
\end{pf}

Let $\mathcal{Q}_L({\zeta_n})$ denote the set of global minimum of
BIC with $\zeta_n$ and its cardinality is $L$.

%
\begin{pro}\label{prou}
Suppose that assumptions \textup{(A1)--(A4)} hold. For $K_n\leq L\leq\widebar{K}_{n}$
and
\[
\Pr\Biggl(\bigcup_{r=1}^{K_n} \bigl\{
\mathcal{Q}_L(\zeta_n)\in D_{r}(L,
\rho_n) \bigr\} \Biggr)\rightarrow0
\]
as $n\rightarrow\infty$, where
\begin{eqnarray*}
&& D_r(L,\rho_n)
\\
&&\qquad = \bigl\{\bigl(\tau_1',
\ldots,\tau_{L}'\bigr)\dvtx 1<\tau_1'<
\cdots<\tau_{L}'\leq n\mbox{ and } \bigl|\tau_s'-
\tau_r\bigr|>\rho_n \mbox{ for } 1\leq s\leq L \bigr\}.
\end{eqnarray*}
\end{pro}

\begin{pf}
For every $(\tau_1',\ldots,\tau_{L}')\in D_r(L,\rho_n)$,
%
%
\begin{eqnarray}
\label{yao11}
&& R_n\bigl(\tau_1',\ldots,
\tau_{L}'\bigr)
\nonumber\\[-8pt]\\[-8pt]
&&\qquad \leq R_n\bigl(
\tau_1',\ldots,\tau_{L}',
\tau_1,\ldots,\tau_{r-1},\tau_r-
\delta_n,\tau_r+\delta_n,
\tau_{r+1},\ldots,\tau_{K_n}\bigr)\nonumber
\end{eqnarray}
and the right-hand side of (\ref{yao11}) can be expressed as
$T_1+\cdots+T_{K_n+2}$, where $T_s$ ($s=1,\ldots,r-1,r+2,\ldots,K_n+1$)
is the sum of squares involving the $X_i$'s ($\tau_{s-1}\leq
i<\tau_s$); $T_r$ is that involving the $X_i$'s ($\tau_{r-1}\leq
i<\tau_r-\rho_n$); $T_{r+1}$ is that involving the $X_i$'s
($\tau_r+\rho_n\leq i<\tau_{r+1}$); $T_{K_n+2}$ is that\vspace*{1pt} involving
the $X_i$'s ($\tau_r-\rho_n\leq i<\tau_r+\rho_n$).
Define $c_n=\widebar{K}_{n}^3(\log\widebar{K}_{n})^2(\log
(n\widebar{K}_{n}))^2$. It can be further seen that uniformly in
$(\tau_1',\ldots,\tau_{L}')\in D_r(L,\rho_n)$,
\begin{eqnarray*}
T_s&=&R_n(\tau_{s-1},\tau_s)+
\widetilde{O}_p(c_n;\widebar{K}_{n}),\qquad s=1,\ldots,r-1,r+2,\ldots,K_n+1,
\\
T_r&=&R_n(\tau_{r-1},\tau_r-
\rho_n)+\widetilde{O}_p(c_n;\widebar
{K}_{n}),
\\
T_{r+1}&=&R_n(\tau_r+\rho_n,
\tau_{r+1})+\widetilde{O}_p(c_n;\widebar
{K}_{n})\quad\mbox{and}
\\
T_{K_n+2}&\leq& R_n(\tau_r-\rho_n,
\tau_r)+R_n(\tau_r,\tau_r+
\rho_n)-\rho_nS(F_{{r}},F_{{r+1}}) +
\widetilde{O}_p(c_n;\widebar{K}_{n}).
\end{eqnarray*}
These results imply that
\[
\mathrm{BIC}_L-\mathrm{BIC}_{*}\geq\rho_nS(F_{{r}},F_{{r+1}})-
\widetilde{O}_p(c_n;\widebar{K}_{n}).
\]
Thus, as $n\to\infty$,
\begin{eqnarray*}
\Pr \Biggl(\bigcup_{r=1}^{K_n} \bigl\{
\mathcal{Q}_L(\zeta_n)\in D_{r}(L,
\rho_n) \bigr\} \Biggr)&\le&\Pr\Biggl(\bigcup
_{r=1}^{K_n}(\mathrm{BIC}_L<
\mathrm{BIC}_{*}) \Biggr)
\\
&\le&\Pr\Biggl(\bigcup_{r=1}^{K_n} \bigl\{
\rho_nS(F_{{r}},F_{{r+1}})<\widetilde{O}_p(c_n;
\widebar{K}_{n}) \bigr\} \Biggr)
\\
&\le& \widebar{K}_{n} \Pr\bigl(\widetilde{O}_p(c_n;
\widebar{K}_{n})>c_n \bigr)<\epsilon
\end{eqnarray*}
for any $\epsilon>0$. Thus, the result follows.
\end{pf}

\begin{pf*}{Proof of Theorem~\ref{teo1}} Define $d_n=K_n^3(\log
K_n)^2(\log(\delta_n K_n))^2$. For every
$(\tau_1',\ldots,\tau_{K_n}')\in D_r(K_n,\delta_n)$,
\begin{eqnarray*}
&& \max_{(\tau_1',\ldots,\tau_{K_n}')\in
D_r(K_n,\delta_n)}R_n\bigl(\tau_1',
\ldots,\tau_{K_n}'\bigr)
\\
&&\qquad \leq R_n\bigl(\tau_1',\ldots,
\tau_{K_n}',\tau_1,\ldots,\tau_{r-1},
\tau_r-\delta_n,\tau_r+
\delta_n,\tau_{r+1},\ldots,\tau_{K_n}\bigr)
\\
&&\qquad \leq R_n(\tau_1,\ldots,\tau_{K_n})-\delta
_nS(F_{{r}},F_{{r+1}})+\widetilde{O}_p(d_n;K_n)
\end{eqnarray*}
by Lemma \ref{klem4}. Thus, we know that
\[
\max_{(\tau_1',\ldots,\tau_{K_n}')\in
D_r(K_n,\delta_n)}R_n\bigl(\tau_1',
\ldots,\tau_{K_n}'\bigr)<R_n(\tau
_1,\ldots,\tau_{K_n})
\]
with probability tending to one for each $r$. Consequently,
\begin{eqnarray*}
\Pr\bigl\{\mathcal{G}_n(K_n)\in C_{K_n}(
\delta_n) \bigr\}&=&1-\Pr\biggl\{\bigcup_r
\bigl\{\mathcal{G}_n(K_n)\in D_r(K_n,
\delta_n) \bigr\} \biggr\}
\\
&\geq&1-\sum_{r=1}^{K_n} \Pr\bigl\{
\mathcal{G}_n(K_n)\in D_r(K_n,
\delta_n) \bigr\} \rightarrow1
\end{eqnarray*}
by the similar argument as that in Proposition \ref{prou}.\vadjust{\goodbreak}
\end{pf*}

\begin{pf*}{Proof of Theorem~\ref{probic}} By Proposition \ref
{prol}, it suffices to show that
$\Pr(\widehat{K}_n>K_n)\rightarrow0$. This can be proved by
contradiction. Let
$E(L,\rho_n)$ be the complement of the union of
$D_1(L,\rho_n),\ldots,D_{K_n}(L,\rho_n)$. As shown in Proposition
\ref{prou}, for $K_n<L<\widebar{K}_{n}$ and every
$(\tau_1',\ldots,\tau_{L}')\in E(L,\rho_n)$,
\begin{eqnarray*}\label{yao1000000}
R_n\bigl(\tau_1',\ldots,
\tau_{L}'\bigr)&\leq& R_n\bigl(
\tau_1',\ldots,\tau_{L}',
\tau_1,\ldots,\tau_{r},\tau_1-
\rho_n, \tau_{K_n}-\rho_n,\tau_1+
\rho_n, \tau_{K_n}+\rho_n\bigr)
\\
&=&R_n(\tau_1,\ldots,\tau_{K_n})+
\widetilde{O}_p(c_n;\widebar{K}_{n}).
\end{eqnarray*}
Consequently, as $n\to\infty$,
\[
\mathrm{BIC}_L- \mathrm{BIC}_{*}\geq(L-K_n)
\zeta_n-\widetilde{O}_p(c_n;
\widebar{K}_{n}),
\]
we obtain the result by the same argument as that in Proposition \ref{prol}.
\end{pf*}
\end{appendix}

\section*{Acknowledgments}
The authors would like to thank Professor Runze Li, an Associate Editor,
and three anonymous
referees for their many insightful and constructive
comments that have resulted in significant improvements in the article.

\begin{supplement}
\stitle{Supplement to ``Nonparametric maximum likelihood approach to
multiple change-point problems''}
\slink[doi]{10.1214/14-AOS1210SUPP} 
\sdatatype{pdf}
\sfilename{AOS1210\_supp.pdf}
\sdescription{We provide technical details for the proof of Corollary~\ref{cor1}, and additional simulation results.}
\end{supplement}

%

\printaddresses

\end{document}